\documentclass{article}

\usepackage{url}
\usepackage{amsmath,amssymb}
\usepackage{color}
\usepackage{tikz}
\usepackage{stmaryrd}
\usepackage{colortbl} 

\usepackage[all]{xypic} \CompileMatrices

\newcommand{\cat}[1]{{\bf #1}}
\newcommand{\comp}{\circ}
\newcommand{\seq}[2]{\shortstack{$#1$ \\ \mbox{}\\
  \mbox{}\hrulefill\mbox{}\\ \mbox{}\\ $#2$}}
\newcommand{\set}[1]{\{#1\}}
\newcommand{\psh}[1]{\cat{Set}^{\cat{#1}^{\it op}}}
\newcommand{\dist}{\mbox{\boldmath $\lambda$}}
\def\profarrow{\!\!\!\xymatrix@C-.75pc{\ar[r]|-{\! +\!} &}\!\!\! }

\newcommand{\Proofitem}[1]{\medskip \noindent $#1\;$}

\newtheorem{theorem}{Theorem}[section]
\newtheorem{definition}[theorem]{Definition}

\newtheorem{proposition}[theorem]{Proposition}

\begin{document}

\title{Operads, clones, and distributive laws}
\author{Pierre-Louis Curien\\ (PPS Laboratory, CNRS \& University Paris Diderot,
Paris, France)}


\maketitle

\begin{abstract}
We show how non-symmetric operads (or multicategories), symmetric operads, and clones, arise from three suitable monads on \cat{Cat}, each extending to a (pseudo-)monad on the bicategory of categories and profunctors. 
We also explain how other previous categorical analyses of operads (via Day's tensor products, or via analytical functors) fit with the profunctor approach. 
\end{abstract}

\section{Introduction}
Operads, in their coloured and non-symmetric version, are also known as multicategories, since they are  like categories, but with morphisms which have a single target as codomain but a multiple source (more precisely, a list of sources) as domain. These morphisms are often called operations (whence the term ``operad'' for the whole structure). Operads for short are the special case where the multicategory has just one object, and hence the arity of a morphism is just a number $n$ (the length of the list), while the coarity is 1. There are many variations.
\begin{itemize}
\item non-symmetric operads versus symmetric operads versus clones (these variations concern the way in which operations are combined to form compound operations);
\item operations versus  cooperations (one input, several outputs), or ``bioperations" (several inputs, several outputs);
\item operations whose shape, or arity, is more structured than a list (it could be a tree, etc...).
\end{itemize}

In this paper, we deal principally with the first variation, and touch on the second briefly. Our goal is
to contribute to convey the idea that these variations can be smoothly and rather uniformly understood using 
some categorical abstractions.  

We are aware of two approaches for a general categorical account of such variations. The first one, which is also the earliest one, is based on spans, while the second one is based on profunctors. In both approaches, one abstracts the details of 
the variant in a monad, and one then extends the monad to the category of spans, or to the category of profunctors. Here, we take the profunctor approach.

\begin{itemize}
\item
Spans are pairs of morphisms with the same domain, in a suitable category \cat{C}. Burroni has shown~[\cite{Burroni-T}]~that every cartesian\footnote{ A monad  $T$ is cartesian if it preserves pullbacks (hence \cat{C} is required to have pullbacks) and if all naturality squares of the unit and multiplication of the monad are pullbacks.}
 monad $T$ on \cat{C} extends to a monad on
spans, and then, for example, multicategories arise as endomorphisms endowed with a monoid structure in the (co)Kleisli (bi)category of a category of spans.
\item
Profunctors are functors $\Phi:\cat{C}\times\cat{C}'^{{\it op}}\rightarrow\cat{Set}$. As pointed out by Cheng~[\cite{Cheng04}]~, every monad in \cat{Cat} (the 2-category of categories) satisfying a certain distributive law can be extended to the category of profunctors (whose objects are categories and whose morphisms are profunctors).  We exhibit here how

\medskip\noindent
{\em non-symmetric and symmetric operads, and clones, arise (again) as monoids 
in the (co)Kleisli (bi)category associated with the respective extended monad.}

\medskip\noindent

\end{itemize}
This is all rather ``heavy" categorical vocabulary. We shall unroll this slowly in what follows. We offer the following intuitions for why monads, spans, and profunctors are relevant here.

\begin{itemize}
\item Monad.
Consider the following composition of operations:
$$
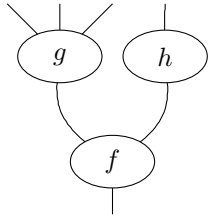
\begin{tikzpicture}[scale=0.7]
\useasboundingbox (-0.5,-0.5) rectangle (3.5,4.5);
\draw (1.00,3.00) -- (0.99,2.95) -- (0.99,2.89) -- (0.98,2.84) -- (0.97,2.79) -- (0.97,2.74) -- (0.96,2.68) -- (0.96,2.63) -- (0.95,2.58) -- (0.95,2.53) -- (0.95,2.48) -- (0.94,2.43) -- (0.94,2.38) -- (0.95,2.33) -- (0.95,2.28) -- (0.95,2.23) -- (0.96,2.18) -- (0.97,2.14) -- (0.97,2.09) -- (0.99,2.04) -- (1.00,2.00);
\draw (1.00,2.00) -- (1.02,1.94) -- (1.05,1.88) -- (1.08,1.82) -- (1.12,1.77) -- (1.15,1.71) -- (1.20,1.66) -- (1.24,1.60) -- (1.29,1.55) -- (1.34,1.50) -- (1.39,1.45) -- (1.44,1.41) -- (1.50,1.36) -- (1.56,1.31) -- (1.62,1.27) -- (1.68,1.22) -- (1.74,1.18) -- (1.81,1.13) -- (1.87,1.09) -- (1.94,1.04) -- (2.00,1.00);
\draw (2.00,4.00) -- (1.00,3.00);
\draw (1.00,4.00) -- (1.00,3.00);
\draw (0.00,4.00) -- (1.00,3.00);
\draw (3.00,4.00) -- (3.00,3.95) -- (3.00,3.90) -- (2.99,3.85) -- (2.99,3.80) -- (2.99,3.75) -- (2.99,3.70) -- (2.99,3.65) -- (2.99,3.61) -- (2.99,3.56) -- (2.99,3.51) -- (2.99,3.46) -- (2.99,3.41) -- (2.99,3.36) -- (2.99,3.31) -- (2.99,3.26) -- (2.99,3.20) -- (2.99,3.15) -- (2.99,3.10) -- (3.00,3.05) -- (3.00,3.00);
\draw (3.00,3.00) -- (3.00,2.95) -- (3.01,2.90) -- (3.01,2.84) -- (3.02,2.79) -- (3.02,2.74) -- (3.03,2.69) -- (3.03,2.64) -- (3.04,2.58) -- (3.04,2.53) -- (3.04,2.48) -- (3.05,2.43) -- (3.05,2.38) -- (3.05,2.33) -- (3.04,2.28) -- (3.04,2.23) -- (3.04,2.18) -- (3.03,2.14) -- (3.02,2.09) -- (3.01,2.04) -- (3.00,2.00);
\draw (3.00,2.00) -- (2.98,1.94) -- (2.95,1.88) -- (2.92,1.82) -- (2.89,1.76) -- (2.85,1.71) -- (2.81,1.65) -- (2.77,1.60) -- (2.72,1.55) -- (2.67,1.50) -- (2.62,1.45) -- (2.56,1.40) -- (2.50,1.36) -- (2.44,1.31) -- (2.38,1.27) -- (2.32,1.22) -- (2.26,1.18) -- (2.20,1.13) -- (2.13,1.09) -- (2.07,1.04) -- (2.00,1.00);
\draw (2.00,1.00) -- (2.00,0.00);
\filldraw[fill=white] (1.00,3.00) ellipse (0.80cm and 0.50cm);
\filldraw[fill=white] (3.00,3.00) ellipse (0.80cm and 0.50cm);
\filldraw[fill=white] (2.00,1.00) ellipse (0.80cm and 0.50cm);
\draw (1.00,3.00) node{$g$};
\draw (3.00,3.00) node{$h$};
\draw (2.00,1.00) node{$f$};
\end{tikzpicture}$$

 By plugging the output of $g$ and $h$ on the inputs of $f$ in parallel, we put together the three inputs of $g$ and the input of $h$, yielding a compound operation with $3+1=4$ inputs. In other words, one may read two shapes on the upper part of the picture:  $((\cdot,\cdot,\cdot),\cdot)$, and $(\cdot,\cdot,\cdot,\cdot)$. The first one
 remembers the construction, the second one flattens it. This flattening is a typical monad multiplication (monad of powersets, of lists, etc...).
  \item Span.  Each bioperation has an arity and a coarity: the legs of the span are this arity and coarity mappings, respectively. 
\item Profunctor: A profunctor gives a family of sets of bioperations of fixed  arity and coarity.
\end{itemize}

Arities (or coarities) are governed by the monad $T$. For example, an operad will be a span / profunctor from $T\cat{1}$ to $\cat{1}$, where $T$ takes care of the multiplicity of inputs (\cat{1} is the category with one object and one morphism). Similarly, a cooperad will go from $1$ to $T\cat{1}$.

The span approach and the profunctor approach should be related, since one goes from profunctors to spans via the ``element" or so-called Grothendieck construction. But under this correspondence the (bi)category of profunctors is (bi)equivalent to a subcategory of spans only, the discrete fibrations, while on the other hand the span approach leaves a lot of freedom on the choice of the underlying category $\cat{C}$.
As a matter of fact, the two approaches have led to different types of successes.
The current state of the art seems to be that:
\begin{enumerate}
\item
using spans, an impressive variety of shapes in the {\em non-symmetric case} have been covered.
We refer to~[\cite{Leinster}]~for a book-length account;
\item using profunctors, one may cover the two other kinds of variations mentioned in this introduction.
\end{enumerate}

In the sequel, we use (and introduce) the profunctor road.
In Section \ref{Burroni-eq}, we introduce several monads on \cat{Cat}. After recalling the notion of Kan extension in Section \ref {Kan-ext}, we present two classical categorical accounts of operads
in Sections \ref{Kelly-op} and \ref{Joyal-op}. We proceed then to 
profunctors in Section \ref{Profunctors}. A  plan of the rest of the paper is given at the end of that section.

\bigskip
\bigskip
\paragraph{Notation}  We shall use juxtaposition to denote functor application.  Moreover, if $F$ is a functor from, say $\cat{C}$ to $\cat{Set}^{\cat{D}}$, we write $FXD$ for $(FX)D$, etc... .
 
\section{Three useful combinators} \label{Burroni-eq}

Consider the following three operations, or combinators. Think of them as morphisms in a category whose objects are (possibly empty) sequences $(C_1,\ldots,C_n)$, where the $C_i$'s range over
the objects of some category \cat{C}.
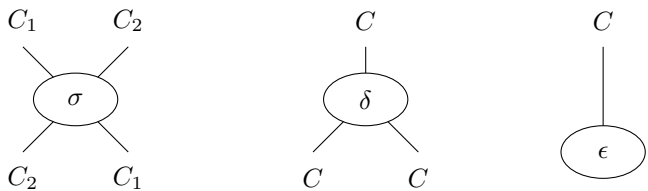
$$\begin{array}{llllllll}
\begin{tikzpicture}[scale=0.7]
\useasboundingbox (-0.5,-0.5) rectangle (2.5,2.5);
\draw (1.00,1.00) -- (2.00,0.00);
\draw (1.00,1.00) -- (0.00,0.00);
\draw (2.00,2.00) -- (1.00,1.00);
\draw (0.00,2.00) -- (1.00,1.00);
\filldraw[fill=white] (1.00,1.00) ellipse (0.80cm and 0.50cm);
\draw (0.00,2.50) node{$C_1$};
\draw (2.00,2.50) node{$C_2$};
\draw (1.00,1.00) node{$\sigma$};
\draw (0.00,-0.50) node{$C_2$};
\draw (2.00,-0.50) node{$C_1$};
\end{tikzpicture}
&\quad\quad&\quad\quad&
\begin{tikzpicture}[scale=0.7]
\useasboundingbox (-0.5,-0.5) rectangle (2.5,2.5);
\draw (1.00,1.00) -- (2.00,0.00);
\draw (1.00,1.00) -- (0.00,0.00);
\draw (1.00,2.00) -- (1.00,1.00);
\filldraw[fill=white] (1.00,1.00) ellipse (0.80cm and 0.50cm);
\draw (1.00,2.50) node{$C$};
\draw (1.00,1.00) node{$\delta$};
\draw (0.00,-0.50) node{$C$};
\draw (2.00,-0.50) node{$C$};
\end{tikzpicture}
&\quad\quad&\quad\quad&
\begin{tikzpicture}[scale=0.7]
\useasboundingbox (-0.5,-0.5) rectangle (0.5,2.5);
\draw (0.00,2.00) -- (0.00,0.00);
\filldraw[fill=white] (0.00,0.00) ellipse (0.80cm and 0.50cm);
\draw (0.00,2.50) node{$C$};
\draw (0.00,0.00) node{$\epsilon$};
\end{tikzpicture}
\end{array}$$

These combinators should satisfy equations:
$$\begin{array}{lllll}
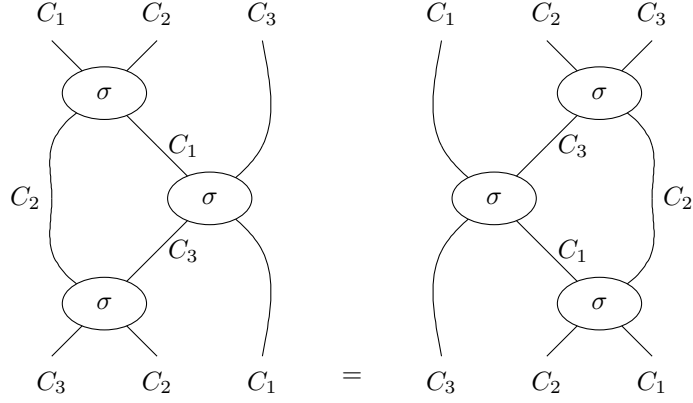
\begin{tikzpicture}[scale=0.7]
\useasboundingbox (-0.5,-0.5) rectangle (4.5,6.5);
\draw (4.00,6.00) -- (4.02,5.89) -- (4.04,5.78) -- (4.06,5.67) -- (4.08,5.57) -- (4.10,5.46) -- (4.11,5.35) -- (4.13,5.25) -- (4.14,5.14) -- (4.15,5.04) -- (4.16,4.94) -- (4.16,4.83) -- (4.16,4.73) -- (4.16,4.64) -- (4.15,4.54) -- (4.14,4.44) -- (4.12,4.35) -- (4.10,4.26) -- (4.07,4.17) -- (4.04,4.08) -- (4.00,4.00);
\draw (4.00,4.00) -- (3.97,3.94) -- (3.94,3.89) -- (3.90,3.83) -- (3.86,3.78) -- (3.82,3.72) -- (3.77,3.67) -- (3.73,3.62) -- (3.68,3.57) -- (3.63,3.52) -- (3.58,3.47) -- (3.52,3.42) -- (3.47,3.37) -- (3.41,3.32) -- (3.36,3.28) -- (3.30,3.23) -- (3.24,3.18) -- (3.18,3.14) -- (3.12,3.09) -- (3.06,3.05) -- (3.00,3.00);
\draw (1.00,5.00) -- (1.05,4.95) -- (1.10,4.90) -- (1.15,4.85) -- (1.20,4.80) -- (1.25,4.75) -- (1.30,4.70) -- (1.35,4.65) -- (1.40,4.60) -- (1.45,4.55) -- (1.50,4.50) -- (1.55,4.45) -- (1.60,4.40) -- (1.65,4.35) -- (1.70,4.30) -- (1.75,4.25) -- (1.80,4.20) -- (1.85,4.15) -- (1.90,4.10) -- (1.95,4.05) -- (2.00,4.00);
\draw (2.00,4.00) -- (2.05,3.95) -- (2.10,3.90) -- (2.15,3.85) -- (2.20,3.80) -- (2.25,3.75) -- (2.30,3.70) -- (2.35,3.65) -- (2.40,3.60) -- (2.45,3.55) -- (2.50,3.50) -- (2.55,3.45) -- (2.60,3.40) -- (2.65,3.35) -- (2.70,3.30) -- (2.75,3.25) -- (2.80,3.20) -- (2.85,3.15) -- (2.90,3.10) -- (2.95,3.05) -- (3.00,3.00);
\draw (1.00,5.00) -- (0.93,4.96) -- (0.87,4.91) -- (0.80,4.87) -- (0.74,4.82) -- (0.67,4.78) -- (0.61,4.73) -- (0.55,4.69) -- (0.49,4.64) -- (0.43,4.59) -- (0.38,4.55) -- (0.32,4.50) -- (0.27,4.45) -- (0.23,4.40) -- (0.18,4.34) -- (0.14,4.29) -- (0.11,4.23) -- (0.07,4.18) -- (0.04,4.12) -- (0.02,4.06) -- (0.00,4.00);
\draw (0.00,4.00) -- (-0.01,3.96) -- (-0.02,3.91) -- (-0.03,3.86) -- (-0.03,3.82) -- (-0.03,3.77) -- (-0.04,3.72) -- (-0.04,3.67) -- (-0.04,3.62) -- (-0.03,3.57) -- (-0.03,3.52) -- (-0.03,3.47) -- (-0.02,3.42) -- (-0.02,3.37) -- (-0.02,3.32) -- (-0.01,3.26) -- (-0.01,3.21) -- (-0.00,3.16) -- (-0.00,3.11) -- (-0.00,3.05) -- (0.00,3.00);
\draw (0.00,3.00) -- (-0.00,2.95) -- (-0.00,2.89) -- (-0.00,2.84) -- (-0.01,2.79) -- (-0.01,2.74) -- (-0.02,2.68) -- (-0.02,2.63) -- (-0.02,2.58) -- (-0.03,2.53) -- (-0.03,2.48) -- (-0.03,2.43) -- (-0.04,2.38) -- (-0.04,2.33) -- (-0.04,2.28) -- (-0.03,2.23) -- (-0.03,2.18) -- (-0.03,2.14) -- (-0.02,2.09) -- (-0.01,2.04) -- (0.00,2.00);
\draw (0.00,2.00) -- (0.02,1.94) -- (0.04,1.88) -- (0.07,1.82) -- (0.11,1.77) -- (0.14,1.71) -- (0.18,1.66) -- (0.23,1.60) -- (0.27,1.55) -- (0.32,1.50) -- (0.38,1.45) -- (0.43,1.41) -- (0.49,1.36) -- (0.55,1.31) -- (0.61,1.27) -- (0.67,1.22) -- (0.74,1.18) -- (0.80,1.13) -- (0.87,1.09) -- (0.93,1.04) -- (1.00,1.00);
\draw (2.00,6.00) -- (1.00,5.00);
\draw (0.00,6.00) -- (1.00,5.00);
\draw (3.00,3.00) -- (3.06,2.95) -- (3.12,2.91) -- (3.18,2.86) -- (3.24,2.82) -- (3.30,2.77) -- (3.36,2.72) -- (3.41,2.68) -- (3.47,2.63) -- (3.52,2.58) -- (3.58,2.53) -- (3.63,2.48) -- (3.68,2.43) -- (3.73,2.38) -- (3.77,2.33) -- (3.82,2.28) -- (3.86,2.22) -- (3.90,2.17) -- (3.94,2.11) -- (3.97,2.06) -- (4.00,2.00);
\draw (4.00,2.00) -- (4.04,1.92) -- (4.07,1.83) -- (4.10,1.74) -- (4.12,1.65) -- (4.14,1.56) -- (4.15,1.46) -- (4.16,1.36) -- (4.16,1.27) -- (4.16,1.17) -- (4.16,1.06) -- (4.15,0.96) -- (4.14,0.86) -- (4.13,0.75) -- (4.11,0.65) -- (4.10,0.54) -- (4.08,0.43) -- (4.06,0.33) -- (4.04,0.22) -- (4.02,0.11) -- (4.00,0.00);
\draw (3.00,3.00) -- (2.95,2.95) -- (2.90,2.90) -- (2.85,2.85) -- (2.80,2.80) -- (2.75,2.75) -- (2.70,2.70) -- (2.65,2.65) -- (2.60,2.60) -- (2.55,2.55) -- (2.50,2.50) -- (2.45,2.45) -- (2.40,2.40) -- (2.35,2.35) -- (2.30,2.30) -- (2.25,2.25) -- (2.20,2.20) -- (2.15,2.15) -- (2.10,2.10) -- (2.05,2.05) -- (2.00,2.00);
\draw (2.00,2.00) -- (1.95,1.95) -- (1.90,1.90) -- (1.85,1.85) -- (1.80,1.80) -- (1.75,1.75) -- (1.70,1.70) -- (1.65,1.65) -- (1.60,1.60) -- (1.55,1.55) -- (1.50,1.50) -- (1.45,1.45) -- (1.40,1.40) -- (1.35,1.35) -- (1.30,1.30) -- (1.25,1.25) -- (1.20,1.20) -- (1.15,1.15) -- (1.10,1.10) -- (1.05,1.05) -- (1.00,1.00);
\draw (1.00,1.00) -- (2.00,0.00);
\draw (1.00,1.00) -- (0.00,0.00);
\filldraw[fill=white] (1.00,5.00) ellipse (0.80cm and 0.50cm);
\filldraw[fill=white] (3.00,3.00) ellipse (0.80cm and 0.50cm);
\filldraw[fill=white] (1.00,1.00) ellipse (0.80cm and 0.50cm);
\draw (0.00,6.50) node{$C_1$};
\draw (2.00,6.50) node{$C_2$};
\draw (4.00,6.50) node{$C_3$};
\draw (1.00,5.00) node{$\sigma$};
\draw (2.50,4.00) node{$C_1$};
\draw (-0.50,3.00) node{$C_2$};
\draw (3.00,3.00) node{$\sigma$};
\draw (2.50,2.00) node{$C_3$};
\draw (1.00,1.00) node{$\sigma$};
\draw (0.00,-0.50) node{$C_3$};
\draw (2.00,-0.50) node{$C_2$};
\draw (4.00,-0.50) node{$C_1$};
\end{tikzpicture}
&& = &&
\begin{tikzpicture}[scale=0.7]
\useasboundingbox (-0.5,-0.5) rectangle (4.5,6.5);
\draw (0.00,6.00) -- (-0.02,5.89) -- (-0.04,5.78) -- (-0.06,5.67) -- (-0.08,5.57) -- (-0.10,5.46) -- (-0.11,5.35) -- (-0.13,5.25) -- (-0.14,5.14) -- (-0.15,5.04) -- (-0.16,4.94) -- (-0.16,4.83) -- (-0.16,4.73) -- (-0.16,4.64) -- (-0.15,4.54) -- (-0.14,4.44) -- (-0.12,4.35) -- (-0.10,4.26) -- (-0.07,4.17) -- (-0.04,4.08) -- (0.00,4.00);
\draw (0.00,4.00) -- (0.03,3.94) -- (0.06,3.89) -- (0.10,3.83) -- (0.14,3.78) -- (0.18,3.72) -- (0.23,3.67) -- (0.27,3.62) -- (0.32,3.57) -- (0.37,3.52) -- (0.42,3.47) -- (0.48,3.42) -- (0.53,3.37) -- (0.59,3.32) -- (0.64,3.28) -- (0.70,3.23) -- (0.76,3.18) -- (0.82,3.14) -- (0.88,3.09) -- (0.94,3.05) -- (1.00,3.00);
\draw (3.00,5.00) -- (3.07,4.96) -- (3.13,4.91) -- (3.20,4.87) -- (3.26,4.82) -- (3.33,4.78) -- (3.39,4.73) -- (3.45,4.69) -- (3.51,4.64) -- (3.57,4.59) -- (3.62,4.55) -- (3.68,4.50) -- (3.73,4.45) -- (3.77,4.40) -- (3.82,4.34) -- (3.86,4.29) -- (3.89,4.23) -- (3.93,4.18) -- (3.96,4.12) -- (3.98,4.06) -- (4.00,4.00);
\draw (4.00,4.00) -- (4.01,3.96) -- (4.02,3.91) -- (4.03,3.86) -- (4.03,3.82) -- (4.03,3.77) -- (4.04,3.72) -- (4.04,3.67) -- (4.04,3.62) -- (4.03,3.57) -- (4.03,3.52) -- (4.03,3.47) -- (4.02,3.42) -- (4.02,3.37) -- (4.02,3.32) -- (4.01,3.26) -- (4.01,3.21) -- (4.00,3.16) -- (4.00,3.11) -- (4.00,3.05) -- (4.00,3.00);
\draw (4.00,3.00) -- (4.00,2.95) -- (4.00,2.89) -- (4.00,2.84) -- (4.01,2.79) -- (4.01,2.74) -- (4.02,2.68) -- (4.02,2.63) -- (4.02,2.58) -- (4.03,2.53) -- (4.03,2.48) -- (4.03,2.43) -- (4.04,2.38) -- (4.04,2.33) -- (4.04,2.28) -- (4.03,2.23) -- (4.03,2.18) -- (4.03,2.14) -- (4.02,2.09) -- (4.01,2.04) -- (4.00,2.00);
\draw (4.00,2.00) -- (3.98,1.94) -- (3.96,1.88) -- (3.93,1.82) -- (3.89,1.77) -- (3.86,1.71) -- (3.82,1.66) -- (3.77,1.60) -- (3.73,1.55) -- (3.68,1.50) -- (3.62,1.45) -- (3.57,1.41) -- (3.51,1.36) -- (3.45,1.31) -- (3.39,1.27) -- (3.33,1.22) -- (3.26,1.18) -- (3.20,1.13) -- (3.13,1.09) -- (3.07,1.04) -- (3.00,1.00);
\draw (3.00,5.00) -- (2.95,4.95) -- (2.90,4.90) -- (2.85,4.85) -- (2.80,4.80) -- (2.75,4.75) -- (2.70,4.70) -- (2.65,4.65) -- (2.60,4.60) -- (2.55,4.55) -- (2.50,4.50) -- (2.45,4.45) -- (2.40,4.40) -- (2.35,4.35) -- (2.30,4.30) -- (2.25,4.25) -- (2.20,4.20) -- (2.15,4.15) -- (2.10,4.10) -- (2.05,4.05) -- (2.00,4.00);
\draw (2.00,4.00) -- (1.95,3.95) -- (1.90,3.90) -- (1.85,3.85) -- (1.80,3.80) -- (1.75,3.75) -- (1.70,3.70) -- (1.65,3.65) -- (1.60,3.60) -- (1.55,3.55) -- (1.50,3.50) -- (1.45,3.45) -- (1.40,3.40) -- (1.35,3.35) -- (1.30,3.30) -- (1.25,3.25) -- (1.20,3.20) -- (1.15,3.15) -- (1.10,3.10) -- (1.05,3.05) -- (1.00,3.00);
\draw (4.00,6.00) -- (3.00,5.00);
\draw (2.00,6.00) -- (3.00,5.00);
\draw (1.00,3.00) -- (1.05,2.95) -- (1.10,2.90) -- (1.15,2.85) -- (1.20,2.80) -- (1.25,2.75) -- (1.30,2.70) -- (1.35,2.65) -- (1.40,2.60) -- (1.45,2.55) -- (1.50,2.50) -- (1.55,2.45) -- (1.60,2.40) -- (1.65,2.35) -- (1.70,2.30) -- (1.75,2.25) -- (1.80,2.20) -- (1.85,2.15) -- (1.90,2.10) -- (1.95,2.05) -- (2.00,2.00);
\draw (2.00,2.00) -- (2.05,1.95) -- (2.10,1.90) -- (2.15,1.85) -- (2.20,1.80) -- (2.25,1.75) -- (2.30,1.70) -- (2.35,1.65) -- (2.40,1.60) -- (2.45,1.55) -- (2.50,1.50) -- (2.55,1.45) -- (2.60,1.40) -- (2.65,1.35) -- (2.70,1.30) -- (2.75,1.25) -- (2.80,1.20) -- (2.85,1.15) -- (2.90,1.10) -- (2.95,1.05) -- (3.00,1.00);
\draw (1.00,3.00) -- (0.94,2.95) -- (0.88,2.91) -- (0.82,2.86) -- (0.76,2.82) -- (0.70,2.77) -- (0.64,2.72) -- (0.59,2.68) -- (0.53,2.63) -- (0.48,2.58) -- (0.42,2.53) -- (0.37,2.48) -- (0.32,2.43) -- (0.27,2.38) -- (0.23,2.33) -- (0.18,2.28) -- (0.14,2.22) -- (0.10,2.17) -- (0.06,2.11) -- (0.03,2.06) -- (0.00,2.00);
\draw (0.00,2.00) -- (-0.04,1.92) -- (-0.07,1.83) -- (-0.10,1.74) -- (-0.12,1.65) -- (-0.14,1.56) -- (-0.15,1.46) -- (-0.16,1.36) -- (-0.16,1.27) -- (-0.16,1.17) -- (-0.16,1.06) -- (-0.15,0.96) -- (-0.14,0.86) -- (-0.13,0.75) -- (-0.11,0.65) -- (-0.10,0.54) -- (-0.08,0.43) -- (-0.06,0.33) -- (-0.04,0.22) -- (-0.02,0.11) -- (0.00,0.00);
\draw (3.00,1.00) -- (4.00,0.00);
\draw (3.00,1.00) -- (2.00,0.00);
\filldraw[fill=white] (3.00,5.00) ellipse (0.80cm and 0.50cm);
\filldraw[fill=white] (1.00,3.00) ellipse (0.80cm and 0.50cm);
\filldraw[fill=white] (3.00,1.00) ellipse (0.80cm and 0.50cm);
\draw (0.00,6.50) node{$C_1$};
\draw (2.00,6.50) node{$C_2$};
\draw (4.00,6.50) node{$C_3$};
\draw (3.00,5.00) node{$\sigma$};
\draw (2.50,4.00) node{$C_3$};
\draw (1.00,3.00) node{$\sigma$};
\draw (4.50,3.00) node{$C_2$};
\draw (2.50,2.00) node{$C_1$};
\draw (3.00,1.00) node{$\sigma$};
\draw (0.00,-0.50) node{$C_3$};
\draw (2.00,-0.50) node{$C_2$};
\draw (4.00,-0.50) node{$C_1$};
\end{tikzpicture}
\end{array}$$
(compare with the familiar equations of transpositions, and notice that here we do not have to index them over natural numbers).
$$\begin{array}{llllllll}
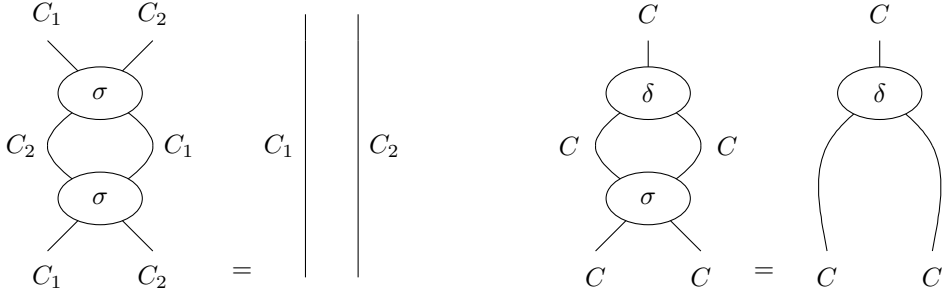
\begin{tikzpicture}[scale=0.7]
\useasboundingbox (-0.5,-0.5) rectangle (2.5,4.5);
\draw (1.00,3.00) -- (1.07,2.95) -- (1.15,2.90) -- (1.22,2.85) -- (1.30,2.80) -- (1.37,2.75) -- (1.44,2.70) -- (1.50,2.65) -- (1.57,2.60) -- (1.63,2.55) -- (1.69,2.50) -- (1.74,2.45) -- (1.79,2.40) -- (1.84,2.35) -- (1.88,2.30) -- (1.91,2.25) -- (1.94,2.20) -- (1.97,2.15) -- (1.99,2.10) -- (2.00,2.05) -- (2.00,2.00);
\draw (2.00,2.00) -- (2.00,1.95) -- (1.99,1.90) -- (1.97,1.85) -- (1.94,1.80) -- (1.91,1.75) -- (1.88,1.70) -- (1.84,1.65) -- (1.79,1.60) -- (1.74,1.55) -- (1.69,1.50) -- (1.63,1.45) -- (1.57,1.40) -- (1.50,1.35) -- (1.44,1.30) -- (1.37,1.25) -- (1.30,1.20) -- (1.22,1.15) -- (1.15,1.10) -- (1.07,1.05) -- (1.00,1.00);
\draw (1.00,3.00) -- (0.93,2.95) -- (0.85,2.90) -- (0.78,2.85) -- (0.70,2.80) -- (0.63,2.75) -- (0.56,2.70) -- (0.50,2.65) -- (0.43,2.60) -- (0.37,2.55) -- (0.31,2.50) -- (0.26,2.45) -- (0.21,2.40) -- (0.16,2.35) -- (0.12,2.30) -- (0.09,2.25) -- (0.06,2.20) -- (0.03,2.15) -- (0.01,2.10) -- (0.00,2.05) -- (0.00,2.00);
\draw (0.00,2.00) -- (0.00,1.95) -- (0.01,1.90) -- (0.03,1.85) -- (0.06,1.80) -- (0.09,1.75) -- (0.12,1.70) -- (0.16,1.65) -- (0.21,1.60) -- (0.26,1.55) -- (0.31,1.50) -- (0.37,1.45) -- (0.43,1.40) -- (0.50,1.35) -- (0.56,1.30) -- (0.63,1.25) -- (0.70,1.20) -- (0.78,1.15) -- (0.85,1.10) -- (0.93,1.05) -- (1.00,1.00);
\draw (2.00,4.00) -- (1.00,3.00);
\draw (0.00,4.00) -- (1.00,3.00);
\draw (1.00,1.00) -- (2.00,0.00);
\draw (1.00,1.00) -- (0.00,0.00);
\filldraw[fill=white] (1.00,3.00) ellipse (0.80cm and 0.50cm);
\filldraw[fill=white] (1.00,1.00) ellipse (0.80cm and 0.50cm);
\draw (0.00,4.50) node{$C_1$};
\draw (2.00,4.50) node{$C_2$};
\draw (1.00,3.00) node{$\sigma$};
\draw (-0.50,2.00) node{$C_2$};
\draw (2.50,2.00) node{$C_1$};
\draw (1.00,1.00) node{$\sigma$};
\draw (0.00,-0.50) node{$C_1$};
\draw (2.00,-0.50) node{$C_2$};
\end{tikzpicture}
&\quad =\quad & 
\begin{tikzpicture}[scale=0.7]
\useasboundingbox (-0.5,-0.5) rectangle (1.5,4.5);
\draw (0.00,4.50) -- (0.00,4.47) -- (0.00,4.45) -- (0.00,4.42) -- (0.00,4.40) -- (0.00,4.38) -- (0.00,4.35) -- (0.00,4.33) -- (0.00,4.30) -- (0.00,4.28) -- (0.00,4.25) -- (0.00,4.22) -- (0.00,4.20) -- (0.00,4.17) -- (0.00,4.15) -- (0.00,4.12) -- (0.00,4.10) -- (0.00,4.08) -- (0.00,4.05) -- (0.00,4.02) -- (0.00,4.00);
\draw (0.00,4.00) -- (0.00,3.77) -- (0.00,3.55) -- (0.00,3.32) -- (0.00,3.10) -- (0.00,2.88) -- (0.00,2.65) -- (0.00,2.42) -- (0.00,2.20) -- (0.00,1.97) -- (0.00,1.75) -- (0.00,1.52) -- (0.00,1.30) -- (0.00,1.07) -- (0.00,0.85) -- (0.00,0.62) -- (0.00,0.40) -- (0.00,0.17) -- (0.00,-0.05) -- (0.00,-0.27) -- (0.00,-0.50);
\draw (1.00,4.50) -- (1.00,4.47) -- (1.00,4.45) -- (1.00,4.42) -- (1.00,4.40) -- (1.00,4.38) -- (1.00,4.35) -- (1.00,4.33) -- (1.00,4.30) -- (1.00,4.28) -- (1.00,4.25) -- (1.00,4.22) -- (1.00,4.20) -- (1.00,4.17) -- (1.00,4.15) -- (1.00,4.12) -- (1.00,4.10) -- (1.00,4.08) -- (1.00,4.05) -- (1.00,4.02) -- (1.00,4.00);
\draw (1.00,4.00) -- (1.00,3.77) -- (1.00,3.55) -- (1.00,3.32) -- (1.00,3.10) -- (1.00,2.88) -- (1.00,2.65) -- (1.00,2.42) -- (1.00,2.20) -- (1.00,1.97) -- (1.00,1.75) -- (1.00,1.52) -- (1.00,1.30) -- (1.00,1.07) -- (1.00,0.85) -- (1.00,0.62) -- (1.00,0.40) -- (1.00,0.17) -- (1.00,-0.05) -- (1.00,-0.27) -- (1.00,-0.50);
\draw (-0.50,2.00) node{$C_1$};
\draw (1.50,2.00) node{$C_2$};
\end{tikzpicture}
&\quad\quad\quad&\quad\quad\quad&
\begin{tikzpicture}[scale=0.7]
\useasboundingbox (-0.5,-0.5) rectangle (2.5,4.5);
\draw (1.00,3.00) -- (1.07,2.95) -- (1.15,2.90) -- (1.22,2.85) -- (1.30,2.80) -- (1.37,2.75) -- (1.44,2.70) -- (1.50,2.65) -- (1.57,2.60) -- (1.63,2.55) -- (1.69,2.50) -- (1.74,2.45) -- (1.79,2.40) -- (1.84,2.35) -- (1.88,2.30) -- (1.91,2.25) -- (1.94,2.20) -- (1.97,2.15) -- (1.99,2.10) -- (2.00,2.05) -- (2.00,2.00);
\draw (2.00,2.00) -- (2.00,1.95) -- (1.99,1.90) -- (1.97,1.85) -- (1.94,1.80) -- (1.91,1.75) -- (1.88,1.70) -- (1.84,1.65) -- (1.79,1.60) -- (1.74,1.55) -- (1.69,1.50) -- (1.63,1.45) -- (1.57,1.40) -- (1.50,1.35) -- (1.44,1.30) -- (1.37,1.25) -- (1.30,1.20) -- (1.22,1.15) -- (1.15,1.10) -- (1.07,1.05) -- (1.00,1.00);
\draw (1.00,3.00) -- (0.93,2.95) -- (0.85,2.90) -- (0.78,2.85) -- (0.70,2.80) -- (0.63,2.75) -- (0.56,2.70) -- (0.50,2.65) -- (0.43,2.60) -- (0.37,2.55) -- (0.31,2.50) -- (0.26,2.45) -- (0.21,2.40) -- (0.16,2.35) -- (0.12,2.30) -- (0.09,2.25) -- (0.06,2.20) -- (0.03,2.15) -- (0.01,2.10) -- (0.00,2.05) -- (0.00,2.00);
\draw (0.00,2.00) -- (0.00,1.95) -- (0.01,1.90) -- (0.03,1.85) -- (0.06,1.80) -- (0.09,1.75) -- (0.12,1.70) -- (0.16,1.65) -- (0.21,1.60) -- (0.26,1.55) -- (0.31,1.50) -- (0.37,1.45) -- (0.43,1.40) -- (0.50,1.35) -- (0.56,1.30) -- (0.63,1.25) -- (0.70,1.20) -- (0.78,1.15) -- (0.85,1.10) -- (0.93,1.05) -- (1.00,1.00);
\draw (1.00,4.00) -- (1.00,3.00);
\draw (1.00,1.00) -- (2.00,0.00);
\draw (1.00,1.00) -- (0.00,0.00);
\filldraw[fill=white] (1.00,3.00) ellipse (0.80cm and 0.50cm);
\filldraw[fill=white] (1.00,1.00) ellipse (0.80cm and 0.50cm);
\draw (1.00,4.50) node{$C$};
\draw (1.00,3.00) node{$\delta$};
\draw (-0.50,2.00) node{$C$};
\draw (2.50,2.00) node{$C$};
\draw (1.00,1.00) node{$\sigma$};
\draw (0.00,-0.50) node{$C$};
\draw (2.00,-0.50) node{$C$};
\end{tikzpicture}
& = & 
\begin{tikzpicture}[scale=0.7]
\useasboundingbox (-0.5,-0.5) rectangle (2.5,4.5);
\draw (1.00,3.00) -- (1.06,2.95) -- (1.12,2.91) -- (1.18,2.86) -- (1.24,2.82) -- (1.30,2.77) -- (1.36,2.72) -- (1.41,2.68) -- (1.47,2.63) -- (1.52,2.58) -- (1.58,2.53) -- (1.63,2.48) -- (1.68,2.43) -- (1.73,2.38) -- (1.77,2.33) -- (1.82,2.28) -- (1.86,2.22) -- (1.90,2.17) -- (1.94,2.11) -- (1.97,2.06) -- (2.00,2.00);
\draw (2.00,2.00) -- (2.04,1.92) -- (2.07,1.83) -- (2.10,1.74) -- (2.12,1.65) -- (2.14,1.56) -- (2.15,1.46) -- (2.16,1.36) -- (2.16,1.27) -- (2.16,1.17) -- (2.16,1.06) -- (2.15,0.96) -- (2.14,0.86) -- (2.13,0.75) -- (2.11,0.65) -- (2.10,0.54) -- (2.08,0.43) -- (2.06,0.33) -- (2.04,0.22) -- (2.02,0.11) -- (2.00,0.00);
\draw (1.00,3.00) -- (0.94,2.95) -- (0.88,2.91) -- (0.82,2.86) -- (0.76,2.82) -- (0.70,2.77) -- (0.64,2.72) -- (0.59,2.68) -- (0.53,2.63) -- (0.48,2.58) -- (0.42,2.53) -- (0.37,2.48) -- (0.32,2.43) -- (0.27,2.38) -- (0.23,2.33) -- (0.18,2.28) -- (0.14,2.22) -- (0.10,2.17) -- (0.06,2.11) -- (0.03,2.06) -- (0.00,2.00);
\draw (0.00,2.00) -- (-0.04,1.92) -- (-0.07,1.83) -- (-0.10,1.74) -- (-0.12,1.65) -- (-0.14,1.56) -- (-0.15,1.46) -- (-0.16,1.36) -- (-0.16,1.27) -- (-0.16,1.17) -- (-0.16,1.06) -- (-0.15,0.96) -- (-0.14,0.86) -- (-0.13,0.75) -- (-0.11,0.65) -- (-0.10,0.54) -- (-0.08,0.43) -- (-0.06,0.33) -- (-0.04,0.22) -- (-0.02,0.11) -- (0.00,0.00);
\draw (1.00,4.00) -- (1.00,3.00);
\filldraw[fill=white] (1.00,3.00) ellipse (0.80cm and 0.50cm);
\draw (1.00,4.50) node{$C$};
\draw (1.00,3.00) node{$\delta$};
\draw (0.00,-0.50) node{$C$};
\draw (2.00,-0.50) node{$C$};
\end{tikzpicture}
\end{array}$$

The other equations are the familiar ones for comonoids, and (mutatis mutandis) for distributive laws (see below).
For a complete list, we refer the reader to~[\cite{Burroni-F,Lafont}]~.

Take any $X\subseteq\set{\sigma,\delta,\epsilon}$. We build a category $!_X\cat{C}$ as follows: objects are sequences $(C_1,\ldots,C_n)$ of objects of \cat{C}~.\footnote{This notation, that comes from linear logic~[\cite{Gir87}]~,
stresses the idea of multiple input (see~[\cite{FGHW}]~for explicit links with linear logic).}
Morphisms are {\em string diagrams} built out of the combinators taken from $X$ and out of the morphisms of  \cat{C}, quotiented over all the equalities that concern the combinators of $X$, including their naturality, and the equalities
$$\begin{array}{llllllllll}
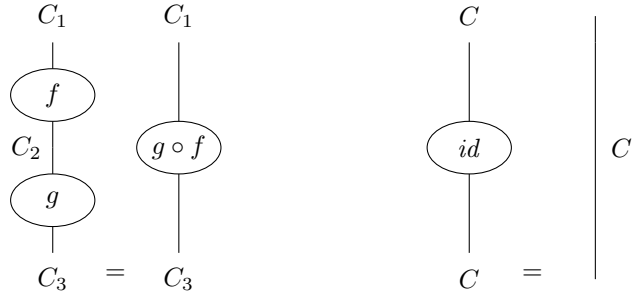
\begin{tikzpicture}[scale=0.7]
\useasboundingbox (-0.5,-0.5) rectangle (0.5,4.5);
\draw (0.00,4.00) -- (0.00,3.95) -- (0.00,3.90) -- (0.00,3.85) -- (0.00,3.80) -- (0.00,3.75) -- (0.00,3.70) -- (0.00,3.65) -- (0.00,3.60) -- (0.00,3.55) -- (0.00,3.50) -- (0.00,3.45) -- (0.00,3.40) -- (0.00,3.35) -- (0.00,3.30) -- (0.00,3.25) -- (0.00,3.20) -- (0.00,3.15) -- (0.00,3.10) -- (0.00,3.05) -- (0.00,3.00);
\draw (0.00,3.00) -- (0.00,2.95) -- (0.00,2.90) -- (0.00,2.85) -- (0.00,2.80) -- (0.00,2.75) -- (0.00,2.70) -- (0.00,2.65) -- (0.00,2.60) -- (0.00,2.55) -- (0.00,2.50) -- (0.00,2.45) -- (0.00,2.40) -- (0.00,2.35) -- (0.00,2.30) -- (0.00,2.25) -- (0.00,2.20) -- (0.00,2.15) -- (0.00,2.10) -- (0.00,2.05) -- (0.00,2.00);
\draw (0.00,2.00) -- (0.00,1.95) -- (0.00,1.90) -- (0.00,1.85) -- (0.00,1.80) -- (0.00,1.75) -- (0.00,1.70) -- (0.00,1.65) -- (0.00,1.60) -- (0.00,1.55) -- (0.00,1.50) -- (0.00,1.45) -- (0.00,1.40) -- (0.00,1.35) -- (0.00,1.30) -- (0.00,1.25) -- (0.00,1.20) -- (0.00,1.15) -- (0.00,1.10) -- (0.00,1.05) -- (0.00,1.00);
\draw (0.00,1.00) -- (0.00,0.95) -- (0.00,0.90) -- (0.00,0.85) -- (0.00,0.80) -- (0.00,0.75) -- (0.00,0.70) -- (0.00,0.65) -- (0.00,0.60) -- (0.00,0.55) -- (0.00,0.50) -- (0.00,0.45) -- (0.00,0.40) -- (0.00,0.35) -- (0.00,0.30) -- (0.00,0.25) -- (0.00,0.20) -- (0.00,0.15) -- (0.00,0.10) -- (0.00,0.05) -- (0.00,0.00);
\filldraw[fill=white] (0.00,3.00) ellipse (0.80cm and 0.50cm);
\filldraw[fill=white] (0.00,1.00) ellipse (0.80cm and 0.50cm);
\draw (0.00,4.50) node{$C_1$};
\draw (0.00,3.00) node{$f$};
\draw (-0.50,2.00) node{$C_2$};
\draw (0.00,1.00) node{$g$};
\draw (0.00,-0.50) node{$C_3$};
\end{tikzpicture}
& = & 
\begin{tikzpicture}[scale=0.7]
\useasboundingbox (-0.5,-0.5) rectangle (0.5,4.5);
\draw (0.00,4.00) -- (0.00,3.90) -- (0.00,3.80) -- (0.00,3.70) -- (0.00,3.60) -- (0.00,3.50) -- (0.00,3.40) -- (0.00,3.30) -- (0.00,3.20) -- (0.00,3.10) -- (0.00,3.00) -- (0.00,2.90) -- (0.00,2.80) -- (0.00,2.70) -- (0.00,2.60) -- (0.00,2.50) -- (0.00,2.40) -- (0.00,2.30) -- (0.00,2.20) -- (0.00,2.10) -- (0.00,2.00);
\draw (0.00,2.00) -- (0.00,1.90) -- (0.00,1.80) -- (0.00,1.70) -- (0.00,1.60) -- (0.00,1.50) -- (0.00,1.40) -- (0.00,1.30) -- (0.00,1.20) -- (0.00,1.10) -- (0.00,1.00) -- (0.00,0.90) -- (0.00,0.80) -- (0.00,0.70) -- (0.00,0.60) -- (0.00,0.50) -- (0.00,0.40) -- (0.00,0.30) -- (0.00,0.20) -- (0.00,0.10) -- (0.00,0.00);
\filldraw[fill=white] (0.00,2.00) ellipse (0.80cm and 0.50cm);
\draw (0.00,4.50) node{$C_1$};
\draw (0.00,2.00) node{$g\comp f$};
\draw (0.00,-0.50) node{$C_3$};
\end{tikzpicture}
&\quad\quad&\quad\quad&\quad\quad&\quad\quad&
\begin{tikzpicture}[scale=0.7]
\useasboundingbox (-0.5,-0.5) rectangle (0.5,4.5);
\draw (0.00,4.00) -- (0.00,3.90) -- (0.00,3.80) -- (0.00,3.70) -- (0.00,3.60) -- (0.00,3.50) -- (0.00,3.40) -- (0.00,3.30) -- (0.00,3.20) -- (0.00,3.10) -- (0.00,3.00) -- (0.00,2.90) -- (0.00,2.80) -- (0.00,2.70) -- (0.00,2.60) -- (0.00,2.50) -- (0.00,2.40) -- (0.00,2.30) -- (0.00,2.20) -- (0.00,2.10) -- (0.00,2.00);
\draw (0.00,2.00) -- (0.00,1.90) -- (0.00,1.80) -- (0.00,1.70) -- (0.00,1.60) -- (0.00,1.50) -- (0.00,1.40) -- (0.00,1.30) -- (0.00,1.20) -- (0.00,1.10) -- (0.00,1.00) -- (0.00,0.90) -- (0.00,0.80) -- (0.00,0.70) -- (0.00,0.60) -- (0.00,0.50) -- (0.00,0.40) -- (0.00,0.30) -- (0.00,0.20) -- (0.00,0.10) -- (0.00,0.00);
\filldraw[fill=white] (0.00,2.00) ellipse (0.80cm and 0.50cm);
\draw (0.00,4.50) node{$C$};
\draw (0.00,2.00) node{${\it id}$};
\draw (0.00,-0.50) node{$C$};
\end{tikzpicture}
& = & 
\begin{tikzpicture}[scale=0.7]
\useasboundingbox (-0.5,-0.5) rectangle (0.5,4.5);
\draw (0.00,4.50) -- (0.00,4.47) -- (0.00,4.45) -- (0.00,4.42) -- (0.00,4.40) -- (0.00,4.38) -- (0.00,4.35) -- (0.00,4.33) -- (0.00,4.30) -- (0.00,4.28) -- (0.00,4.25) -- (0.00,4.22) -- (0.00,4.20) -- (0.00,4.17) -- (0.00,4.15) -- (0.00,4.12) -- (0.00,4.10) -- (0.00,4.08) -- (0.00,4.05) -- (0.00,4.02) -- (0.00,4.00);
\draw (0.00,4.00) -- (0.00,3.77) -- (0.00,3.55) -- (0.00,3.32) -- (0.00,3.10) -- (0.00,2.88) -- (0.00,2.65) -- (0.00,2.42) -- (0.00,2.20) -- (0.00,1.97) -- (0.00,1.75) -- (0.00,1.52) -- (0.00,1.30) -- (0.00,1.07) -- (0.00,0.85) -- (0.00,0.62) -- (0.00,0.40) -- (0.00,0.17) -- (0.00,-0.05) -- (0.00,-0.27) -- (0.00,-0.50);
\draw (0.50,2.00) node{$C$};
\end{tikzpicture}
\end{array}$$
 which allow to embed $\cat{C}$ functorially into $!_X\cat{C}$. String diagrams are combinations of vertical and horizontal combinations of the basic combinators $\sigma,\delta,\epsilon,f$, such as the ones drawn above.
 
 \medskip\noindent
 
 For six of the eight choices of $X$,
 the functor $!_X$ is equipped with a monad structure: the unit takes $C$ to $(C)$ and $f:C\rightarrow D$ to $f:(C)\rightarrow(D)$ and the multiplication on objects is the usual flattening, which takes
 $((C_1^1,\ldots,C_1^{i_1}),(C_n^1,\ldots,C_n^{i_n}))$ to $(C_1^1,\ldots,C_n^{i_n})$.
 A little care is needed to define the ``flattening" on morphisms. A morphism of $!!\cat{C}$ is an assembling of boxes connected by combinators typed in
 $!\cat{C}$.  When we remove the boxes, we need to turn these combinators into (assemblings of) combinators typed in \cat{C}. We call this an expansion.  A $\sigma$ is expanded by means of $\sigma$'s, an $\epsilon$ is expanded by putting $\epsilon$'s in parallel, but we need $\delta$'s {\em and} $\sigma$'s to expand a $\delta$. For example, at type $(C_1,C_2)$, the expansion of 
 $\delta:((C_1,C_2))\rightarrow ((C_1,C_2),(C_1,C_2))$ is the following morphism from
 $(C_1,C_2)$ to $(C_1,C_2,C_1,C_2)$:
$$
\begin{tikzpicture}[scale=0.7]
\useasboundingbox (-0.5,-0.5) rectangle (6.5,4.5);
\draw (1.00,3.00) -- (1.05,2.95) -- (1.10,2.90) -- (1.15,2.85) -- (1.20,2.80) -- (1.25,2.75) -- (1.30,2.70) -- (1.35,2.65) -- (1.40,2.60) -- (1.45,2.55) -- (1.50,2.50) -- (1.55,2.45) -- (1.60,2.40) -- (1.65,2.35) -- (1.70,2.30) -- (1.75,2.25) -- (1.80,2.20) -- (1.85,2.15) -- (1.90,2.10) -- (1.95,2.05) -- (2.00,2.00);
\draw (2.00,2.00) -- (2.05,1.95) -- (2.10,1.90) -- (2.15,1.85) -- (2.20,1.80) -- (2.25,1.75) -- (2.30,1.70) -- (2.35,1.65) -- (2.40,1.60) -- (2.45,1.55) -- (2.50,1.50) -- (2.55,1.45) -- (2.60,1.40) -- (2.65,1.35) -- (2.70,1.30) -- (2.75,1.25) -- (2.80,1.20) -- (2.85,1.15) -- (2.90,1.10) -- (2.95,1.05) -- (3.00,1.00);
\draw (1.00,3.00) -- (0.94,2.95) -- (0.88,2.91) -- (0.82,2.86) -- (0.76,2.82) -- (0.70,2.77) -- (0.64,2.72) -- (0.59,2.68) -- (0.53,2.63) -- (0.48,2.58) -- (0.42,2.53) -- (0.37,2.48) -- (0.32,2.43) -- (0.27,2.38) -- (0.23,2.33) -- (0.18,2.28) -- (0.14,2.22) -- (0.10,2.17) -- (0.06,2.11) -- (0.03,2.06) -- (0.00,2.00);
\draw (0.00,2.00) -- (-0.04,1.92) -- (-0.07,1.83) -- (-0.10,1.74) -- (-0.12,1.65) -- (-0.14,1.56) -- (-0.15,1.46) -- (-0.16,1.36) -- (-0.16,1.27) -- (-0.16,1.17) -- (-0.16,1.06) -- (-0.15,0.96) -- (-0.14,0.86) -- (-0.13,0.75) -- (-0.11,0.65) -- (-0.10,0.54) -- (-0.08,0.43) -- (-0.06,0.33) -- (-0.04,0.22) -- (-0.02,0.11) -- (0.00,0.00);
\draw (1.00,4.00) -- (1.00,3.00);
\draw (5.00,3.00) -- (5.06,2.95) -- (5.12,2.91) -- (5.18,2.86) -- (5.24,2.82) -- (5.30,2.77) -- (5.36,2.72) -- (5.41,2.68) -- (5.47,2.63) -- (5.52,2.58) -- (5.58,2.53) -- (5.63,2.48) -- (5.68,2.43) -- (5.73,2.38) -- (5.77,2.33) -- (5.82,2.28) -- (5.86,2.22) -- (5.90,2.17) -- (5.94,2.11) -- (5.97,2.06) -- (6.00,2.00);
\draw (6.00,2.00) -- (6.04,1.92) -- (6.07,1.83) -- (6.10,1.74) -- (6.12,1.65) -- (6.14,1.56) -- (6.15,1.46) -- (6.16,1.36) -- (6.16,1.27) -- (6.16,1.17) -- (6.16,1.06) -- (6.15,0.96) -- (6.14,0.86) -- (6.13,0.75) -- (6.11,0.65) -- (6.10,0.54) -- (6.08,0.43) -- (6.06,0.33) -- (6.04,0.22) -- (6.02,0.11) -- (6.00,0.00);
\draw (5.00,3.00) -- (4.95,2.95) -- (4.90,2.90) -- (4.85,2.85) -- (4.80,2.80) -- (4.75,2.75) -- (4.70,2.70) -- (4.65,2.65) -- (4.60,2.60) -- (4.55,2.55) -- (4.50,2.50) -- (4.45,2.45) -- (4.40,2.40) -- (4.35,2.35) -- (4.30,2.30) -- (4.25,2.25) -- (4.20,2.20) -- (4.15,2.15) -- (4.10,2.10) -- (4.05,2.05) -- (4.00,2.00);
\draw (4.00,2.00) -- (3.95,1.95) -- (3.90,1.90) -- (3.85,1.85) -- (3.80,1.80) -- (3.75,1.75) -- (3.70,1.70) -- (3.65,1.65) -- (3.60,1.60) -- (3.55,1.55) -- (3.50,1.50) -- (3.45,1.45) -- (3.40,1.40) -- (3.35,1.35) -- (3.30,1.30) -- (3.25,1.25) -- (3.20,1.20) -- (3.15,1.15) -- (3.10,1.10) -- (3.05,1.05) -- (3.00,1.00);
\draw (5.00,4.00) -- (5.00,3.00);
\draw (3.00,1.00) -- (4.00,0.00);
\draw (3.00,1.00) -- (2.00,0.00);
\filldraw[fill=white] (1.00,3.00) ellipse (0.80cm and 0.50cm);
\filldraw[fill=white] (5.00,3.00) ellipse (0.80cm and 0.50cm);
\filldraw[fill=white] (3.00,1.00) ellipse (0.80cm and 0.50cm);
\draw (1.00,4.50) node{$C_1$};
\draw (5.00,4.50) node{$C_2$};
\draw (1.00,3.00) node{$\delta$};
\draw (5.00,3.00) node{$\delta$};
\draw (1.50,2.00) node{$C_1$};
\draw (4.50,2.00) node{$C_2$};
\draw (3.00,1.00) node{$\sigma$};
\draw (0.00,-0.50) node{$C_1$};
\draw (2.00,-0.50) node{$C_2$};
\draw (4.00,-0.50) node{$C_1$};
\draw (6.00,-0.50) node{$C_2$};
\end{tikzpicture}
$$
We therefore exclude $!_{\set{\delta}}$ and $!_{\set{\delta,\epsilon}}$ from our treatment. But all the other $!$ constructions are fine. 

 Three of the six remaining combinations are of particular interest.
 If  $X=\emptyset$ (resp. $X=\set{\sigma}$, $X=\set{\sigma,\delta,\epsilon}$), then $!\cat{C}$ is the free strict  monoidal (resp. symmetric monoidal, cartesian) category over $\cat{C}$. (By cartesian category, we mean ``category with specified finite products").
We write them $!_m,!_s,!_f$, respectively:
\begin{center}
\begin{tabular}{ccc}
monoidal & symmetric monoidal & cartesian\\\\
$!_m$ & $!_s$ & $!_f$
\end{tabular}
\end{center}
When  $X=\set{\sigma,\delta,\epsilon}$, the $!$-algebras are the cartesian categories, and when $X=\emptyset$ (resp. $X=\set{\sigma}$), the  (pseudo-) $!$-algebras are the monoidal (resp. symmetric monoidal) categories. We refer to~[\cite{Leinster}]~for details.

\medskip
We define $?\cat{C}=(!(\cat{C}^{\it op}))^{\it op}$.
Graphically, this amounts to reversing the basic combinators $\sigma,\delta,\epsilon$ (which could then, if we cared, be called $\sigma,\mu,\eta$), while maintaining the orientation of the combinators $f$ imported from $\cat{C}$ (since they become themselves again after two ${}^{\it op}$'s). But it is more convenient to stick with $\sigma,\delta,\epsilon$, to reverse the direction of the $f$'s, and to read the diagrams in the bottom-to-top direction. In particular, when $\cat{C}=\cat{1}$ (the terminal category), there is no $f$ to reverse...

When $X= \emptyset$ or  $X=\set{\sigma}$, then $?\cat{C}$ is (isomorphic to) $!\cat{C}$.
When $X=\set{\sigma,\delta,\epsilon}$, $?\cat{C}$ is the free cocartesian category over \cat{C}.
The objects of $?\cat{1}$ are written $0,1,2,\ldots$, standing for $(),(\cdot),(\cdot,\cdot,\cdot),\ldots$, where $\cdot$ is the unique object of $\cat{1}$.

Note that for any choice of $?$, there is a faithful functor from $?\cat{1}$ to \cat{Set}.
For any morphism from $m$ to $n$, i.e., for any string diagram
constructed out of the combinators in $X$ that has $n$ input wires and $m$ output wires, one constructs a function from $\set{0,\ldots,m-1}$ to $\set{0,\ldots,n-1}$, by the following rules.
\begin{itemize}
\item $\sigma$: the transposition on $\set{0,1}$;
\item $\epsilon$:  the unique function from the empty set to $\set{0}$;
\item $\delta$: the unique function from $\set{0,1}$ to $\set{0}$;
\item vertical composition of string diagrams: function composition;
\item horizontal composition of string diagrams: their categorical sum, i.e., for $f:\set{0,\ldots,m-1}\rightarrow\set{0,\ldots,n-1}$ and $g:\set{0,\ldots,p-1}\rightarrow\set{0,\ldots,q-1}$,  $f+g:\set{0,\ldots,m+p-1}\rightarrow\set{0,\ldots,n+q-1}$ is defined by  $(f+g)(i)=f(i)$ if $i<m$ and $f+g(i)=g(i)+n$ if $i\geq m$.
\end{itemize}

More synthetically, the function $f$ associated with a string diagram is obtained by naming the output and input  wires 
$0,1,\ldots,m-1$ and $0,1,\ldots,n-1$, respectively, and then computing $f(i)$ as the name of the input wire reached when starting from input wire $i$, going up in $\delta$ nodes, and  going up left or right in $\sigma$ nodes according to whether this node is reached from down right or down left, respectively.
For example, the picture
$$
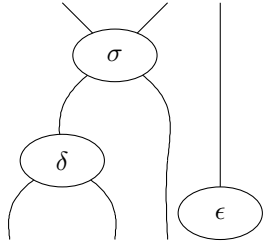
\begin{tikzpicture}[scale=0.7]
\useasboundingbox (-0.5,-0.5) rectangle (4.5,4.5);
\draw (2.00,3.00) -- (2.07,2.96) -- (2.13,2.91) -- (2.19,2.87) -- (2.26,2.82) -- (2.32,2.78) -- (2.38,2.73) -- (2.44,2.69) -- (2.50,2.64) -- (2.56,2.60) -- (2.61,2.55) -- (2.67,2.50) -- (2.72,2.45) -- (2.76,2.40) -- (2.81,2.35) -- (2.85,2.29) -- (2.89,2.24) -- (2.92,2.18) -- (2.95,2.12) -- (2.98,2.06) -- (3.00,2.00);
\draw (3.00,2.00) -- (3.01,1.96) -- (3.02,1.91) -- (3.03,1.86) -- (3.04,1.82) -- (3.04,1.77) -- (3.05,1.72) -- (3.05,1.67) -- (3.05,1.62) -- (3.05,1.57) -- (3.05,1.52) -- (3.04,1.47) -- (3.04,1.42) -- (3.04,1.36) -- (3.03,1.31) -- (3.03,1.26) -- (3.02,1.21) -- (3.02,1.16) -- (3.01,1.10) -- (3.00,1.05) -- (3.00,1.00);
\draw (3.00,1.00) -- (2.99,0.92) -- (2.99,0.85) -- (2.98,0.77) -- (2.98,0.69) -- (2.98,0.62) -- (2.98,0.54) -- (2.97,0.46) -- (2.97,0.39) -- (2.97,0.31) -- (2.97,0.24) -- (2.98,0.16) -- (2.98,0.09) -- (2.98,0.02) -- (2.98,-0.06) -- (2.98,-0.13) -- (2.99,-0.21) -- (2.99,-0.28) -- (2.99,-0.35) -- (3.00,-0.43) -- (3.00,-0.50);
\draw (2.00,3.00) -- (1.94,2.96) -- (1.87,2.91) -- (1.81,2.87) -- (1.74,2.82) -- (1.68,2.78) -- (1.62,2.73) -- (1.56,2.69) -- (1.50,2.64) -- (1.44,2.59) -- (1.39,2.55) -- (1.34,2.50) -- (1.29,2.45) -- (1.24,2.40) -- (1.20,2.34) -- (1.15,2.29) -- (1.12,2.23) -- (1.08,2.18) -- (1.05,2.12) -- (1.02,2.06) -- (1.00,2.00);
\draw (1.00,2.00) -- (0.99,1.96) -- (0.97,1.91) -- (0.97,1.86) -- (0.96,1.82) -- (0.95,1.77) -- (0.95,1.72) -- (0.95,1.67) -- (0.94,1.62) -- (0.94,1.57) -- (0.95,1.52) -- (0.95,1.47) -- (0.95,1.42) -- (0.96,1.37) -- (0.96,1.32) -- (0.97,1.26) -- (0.97,1.21) -- (0.98,1.16) -- (0.99,1.11) -- (0.99,1.05) -- (1.00,1.00);
\draw (3.00,4.00) -- (2.00,3.00);
\draw (1.00,4.00) -- (2.00,3.00);
\draw (1.00,1.00) -- (1.07,0.96) -- (1.14,0.92) -- (1.20,0.87) -- (1.27,0.83) -- (1.34,0.79) -- (1.40,0.74) -- (1.46,0.70) -- (1.52,0.65) -- (1.58,0.60) -- (1.64,0.56) -- (1.69,0.51) -- (1.74,0.46) -- (1.79,0.41) -- (1.83,0.35) -- (1.87,0.30) -- (1.91,0.24) -- (1.94,0.19) -- (1.96,0.13) -- (1.98,0.06) -- (2.00,0.00);
\draw (2.00,0.00) -- (2.00,-0.02) -- (2.01,-0.05) -- (2.01,-0.07) -- (2.01,-0.09) -- (2.02,-0.12) -- (2.02,-0.14) -- (2.02,-0.17) -- (2.02,-0.19) -- (2.02,-0.22) -- (2.02,-0.24) -- (2.02,-0.27) -- (2.02,-0.29) -- (2.01,-0.32) -- (2.01,-0.34) -- (2.01,-0.37) -- (2.01,-0.40) -- (2.01,-0.42) -- (2.00,-0.45) -- (2.00,-0.47) -- (2.00,-0.50);
\draw (1.00,1.00) -- (0.93,0.96) -- (0.86,0.92) -- (0.80,0.87) -- (0.73,0.83) -- (0.66,0.79) -- (0.60,0.74) -- (0.54,0.70) -- (0.48,0.65) -- (0.42,0.60) -- (0.36,0.56) -- (0.31,0.51) -- (0.26,0.46) -- (0.21,0.41) -- (0.17,0.35) -- (0.13,0.30) -- (0.09,0.24) -- (0.06,0.19) -- (0.04,0.13) -- (0.02,0.06) -- (0.00,0.00);
\draw (0.00,0.00) -- (-0.00,-0.02) -- (-0.01,-0.05) -- (-0.01,-0.07) -- (-0.01,-0.09) -- (-0.02,-0.12) -- (-0.02,-0.14) -- (-0.02,-0.17) -- (-0.02,-0.19) -- (-0.02,-0.22) -- (-0.02,-0.24) -- (-0.02,-0.27) -- (-0.02,-0.29) -- (-0.01,-0.32) -- (-0.01,-0.34) -- (-0.01,-0.37) -- (-0.01,-0.40) -- (-0.01,-0.42) -- (-0.00,-0.45) -- (-0.00,-0.47) -- (0.00,-0.50);
\draw (4.00,4.00) -- (4.00,0.00);
\filldraw[fill=white] (2.00,3.00) ellipse (0.80cm and 0.50cm);
\filldraw[fill=white] (1.00,1.00) ellipse (0.80cm and 0.50cm);
\filldraw[fill=white] (4.00,0.00) ellipse (0.80cm and 0.50cm);
\draw (2.00,3.00) node{$\sigma$};
\draw (1.00,1.00) node{$\delta$};
\draw (4.00,0.00) node{$\epsilon$};
\end{tikzpicture}
$$
represents the function $f:\set{0,1,2}\rightarrow\set{0,1,2}$ defined by
$$f(0)=f(1)=1\quad f(2)=0$$
Note on this example that $\epsilon$ witnesses the lack of surjectivity.
More precisely, as first shown by Burroni~[\cite{Burroni-F}]~, the functor, when suitably corestricted, gives the following isomorphims and  equivalences of categories:
$$\begin{array}{lll}
?_\emptyset && \mbox{identity functions}\\
?_{\set{\sigma}} && \mbox{bijections}\\
?_{\set{\sigma,\delta,\epsilon}} && \mbox{all functions}\\
?_{\set{\delta,\epsilon}} && \mbox{monotone functions}\\
?_{\set{\sigma,\delta}} && \mbox{surjective functions}\\
?_{\set{\sigma,\epsilon}} && \mbox{injective functions}
\end{array}$$
where, say, the last line should be read as follows: the category $?_{\set{\sigma,\epsilon}}(\cat{1})$ is equivalent  to the category of finite sets and injective functions, and isomorphic to its full subcategory whose objects are the sets $\set{0,\ldots,i-1}$ ($i\in\mathbb{N}$). For $?_{\set{\delta,\epsilon}}(\cat{1})$, one must take the category of all finite total orders.

 \section{Kan extensions} \label{Kan-ext}
 
 We recall that given a functor $K:\cat{M}\rightarrow\cat{C}$, a left Kan extension of a functor
 $T:\cat{M}\rightarrow\cat{A}$ along $K$ is a pair of a functor ${\it Lan}_K(T): \cat{C}\rightarrow\cat{A}$ and a natural  transformation $\eta_T:T\rightarrow {\it Lan}_K(T)K$ such that $({\it Lan}_K(T),\eta_T)$ is  universal from  $T$ to $\cat{A}^K$. This means that for any other pair of a functor 
 $S:\cat{C}\rightarrow\cat{A}$ and a natural transformation $\alpha:T\rightarrow SK$ there exists a unique 
natural transformation $\mu:{\it Lan}_K(T)\rightarrow S$ such that $\alpha=\mu K\comp \eta$. Here we shall only need to know that the solution of this universal problem is given by the following formula for $ {\it Lan}_K(T)$:
$${\it Lan}_KTC=\int^m\cat{C}[KM,C]\cdot TM$$
where we use Mac Lane's notation for coends~[\cite{MacLane}]~. Coends are sorts of colimits (or inductive limits), adapted to the case of diagrams which vary both covariantly and contravariantly over some parameter: here, $M$ appears contravariantly in $\cat{C}[KM,C]$ and covariantly in $TM$.
In the formula, $\cat{C}[KM,C]\cdot TM$ stands for the coproduct of as many copies of $TM$ as there are 
morphisms from $KM$ to $C$.
 
$$  \xymatrix @C=0.4cm  {\cat{M} \ar[d]^{T} \ar[r]^K &  \cat{C}\ar[dl]^{{\it Lan}_KT}\\
\cat{A} &} 
$$

When $K$ is full and faithful, then $\eta$ is iso, and in particular  the triangle commutes (up to isomorphism).
When moreover $K$ is the Yoneda embedding, we get
$${\it Lan}_{\cal Y}TX=\int^MXm\cdot TM$$
When moreover $\cat{A}$ is a presheaf category $\cat{Set}^\cat{D}$, we have (since limits and colimits of presheaves are pointwise):
$${\it Lan}_{\cal Y}TXD=\int^MXM\times TMD$$
or, making the quotient involved in the coend explicit:
$${\it Lan}_{\cal Y}TXD=(\sum_MXM\times TMD)/\approx$$
 When $\approx$ is termwise, i.e. is the disjoint union of equivalences $\approx_M$ each on 
$XM\times TMD$,  we get a  formula which will look more familiar to algebraic operadists:
$${\it Lan}_{\cal Y}TXD=\sum_M (XM\otimes_M TMD)$$
(where we have written $(XM\times TMD)/{\approx_M}$ as $XM\otimes_MTMD$).
In our examples, $\approx$ will not always be termwise (it will be in the operad cases, but not in the
clone case).

\section{Kelly's account of operads} \label{Kelly-op}

In 1972, Kelly gave the following description of operads~[\cite{Kelly-Operads}]~(see also~[\cite{FioreFossacs}]~):

$$
 \xymatrix @C=0.4cm  {\cat{1} \ar[drr]_{X} \ar[rr] && !\cat{1} \ar[d]^{X^{\otimes\_}} \ar[rr]^{\cal{Y}} && \cat{Set}^{?\cat{1}}\ar[dll]^{\;\;\;{\it Lan}_{\cal{Y}}(X^{\otimes\_})=\_\bullet X}\\
 && \cat{Set}^{?\cat{1}} && } 
 $$
In this diagram, the functor $X^{\otimes \_}$ associates with $n$ the $n$th iterated tensor of $X$, with respect to the following tensor product  structure on $\cat{Set}^{?\cat{1}}$ due (in a more general setting) to Day~[\cite{Day}]~:
$$(X\otimes Y)p=\int^{m,n} Xm\times Yn\times  ?\cat{1}[m+n,p]
$$
The $n$-ary tensor product is described by the following formula (which also gives the unit of the tensor product, for $m=0$):
$$(X_1\otimes X_2\otimes\ldots\otimes X_m)p=\int^{n_1,\ldots,n_m} X_1n_1\times\ldots\times X_mn_m\times  ?\cat{1}[n_1+\ldots+n_m,p]$$
Hence the formula for ${\it Lan}_{\cal{Y}}(X^{\otimes\_})Y$, which we write $Y\bullet X$, is
$$\begin{array}{lll}
(Y\bullet X)p  & = & \int^m Ym\times X^{\otimes m}p\\ 
& = & \int^{m,n_1,\ldots,n_m}Ym\times Xn_1\times\ldots\times Xn_m\times?\cat{1}[n_1+\ldots+n_m,p]
\end{array}$$
In the case of $!_f$ this boils down to
$$(Y\bullet X)p=\int^mYm\times (Xp)^m$$
This is because in a cocartesian category, a morphism $f:n_1+\ldots+n_m\rightarrow p$ amounts to morphisms $f_i\in?\cat{1}[n_i,p]$ which allows us to define  a map from  $Xn_1\times\ldots\times Xn_m\times?\cat{1}[n_1+\ldots+n_m,p]$ to $(Xp)^m$ as follows:
$$(a_1,\ldots,a_m,f)\mapsto (Xf_1a_1,\ldots,Xf_ma_m)$$
and then to ``extend'' cocones indexed over $m$ to cocones indexed over $m,n_1,\ldots,n_m$.

\smallskip
One can check rather easily that the {\em substitutions operation} $\bullet$, together with the unit defined by
$$In=?\cat{1}[1,n]$$
provide a (non-symmetric) monoidal structure on $\cat{Set}^{?\cat{1}}$
(for $I\bullet X\approx X$, one uses the fact that $X^{\otimes\_}$ is functorial).

When ! is $!_m,!_s,!_f$, respectively, a monoid for this structure $(I,\bullet)$ is an operad, a symmetric operad, a
clone, respectively.

\section{Operads from analytic functors} \label{Joyal-op}

Another way to arrive at  the operation $\bullet$ just defined is via a different Kan extension:
$$
 \xymatrix @C=0.4cm  {\cat{?\cat{1}} \ar[r]^{\subseteq} \ar[d]^X &  \cat{Set}\ar[dl]^{{\it Lan}_\subseteq X}\\
\cat{Set} &} 
$$
where $\subseteq $ is the faithful (and non full) functor described at the end of Section \ref{Burroni-eq}.
This is the approach taken by Joyal (for $!_s$)~[\cite{Joyal2}]~. In Joyal's language, $X$ is called a {\em species of structure}, and ${\it Lan}_\subseteq X$ is the associated {\em analytic functor}, whose explicit formula is (for any set $z$):
$${\it Lan}_\subseteq Xz=\int^m z^m\times Xm$$
It can be shown  that ${\it Lan}_\subseteq$ is faithful.  Then $Y\bullet X$ is characterised by the following equality:
$${\it Lan}_\subseteq(Y\bullet X)={\it Lan}_\subseteq Y\comp{\it Lan}_\subseteq X$$
which evidences the fact that $\bullet$ is a composition operation.
Indeed, we have:
$$\begin{array}{l}
{\it Lan}_\subseteq Y({\it Lan}_\subseteq Xz)\\
\quad = \int^m(\int^n z^n\times Xn)^m\times Ym\\
\quad = \int^{m,n_1,\ldots,n_m}z^{n_1+\ldots+n_m}\times Xn_1\times\ldots\times Xn_m\times Ym\\
\quad = \int^p z^p\times(\int^{m,n_1,\ldots,n_m} Ym\times Xn_1\times\ldots\times Xn_m\times ?\cat{1}[n_1+\ldots+n_m,p])\\
\quad = {\it Lan}_\subseteq(Y\bullet X)z
\end{array}$$
(the summation over $p$ is superfluous since from $f\in?\cat{1}[n_1+\ldots+n_m,p]$ we get $z^{\subseteq(f)}$ from $z^p$ to $z^{n_1+\ldots+n_m}$).

 \section{Profunctors} \label{Profunctors}
 Recall that a profunctor (or distributor)~[\cite{BenabouDis,Borceux1}]~$\Phi$ from $\cat{C}$ to $\cat{C'}$ is a functor
$$\Phi:\cat{C}\times\cat{C'}^{\it op}\rightarrow \cat{Set}$$ 
 We write $\Phi:\cat{C}\profarrow \cat{C'}$.
Composition of profunctors is given by the following formula:
$$(\Psi\comp\Phi)(C,C'')=\int^{C'}\Psi(C',C'')\times\Phi(C,C')$$
Therefore, profunctors compose only up to isomorphism. Categories, profunctors, and natural transformations form thus, not a 2-category, but a bicategory~[\cite{BenabouBi}]~.

The bicategory \cat{Prof} of profunctors is self-dual, via the isomorphism ${}^{\it op}:\cat{Prof}\rightarrow\cat{Prof}^{\it op}$ which maps $\cat{C}$ to $\cat{C}^{\it op}$ and
$\Phi:\cat{C_1}\profarrow\cat{C_2}$ to 
$$\Phi^{\it op}=((C_2,C_1)\mapsto \Phi(C_1,C_2)): \cat{C_2}^{\it op}\profarrow\cat{C_1}^{\it op}$$

The composition of profunctors can be synthesised via Kan extensions: Given $\Phi:\cat{C}\profarrow \cat{C'}$ and $\Psi:\cat{C'}\profarrow\cat{C''}$, consider their ``twisted curried"\footnote{After the name of Curry,  who defined a calculus of functions called combinatory logic, based on application and a few combinators, and where functions of several arguments are expressed through repeated applications.}
 versions
$\Phi':\cat{C'}^{\it op}\rightarrow\cat{Set}^\cat{C}$ and
$\Psi':\cat{C''}^{\it op}\rightarrow\cat{Set}^\cat{C'}$, defined by $\Phi'C'C=\Phi(C,C')$ and
$\Psi'C''C'=\Psi(C',C'')$. Then it is immediate to check that
$\Psi\comp\Phi$ is the uncurried version of
$({\it Lan}_{\cal Y}\Phi')\comp\Psi'$, as illustrated in the following diagram:
$$
\xymatrix @C=0.4cm  {\cat{C'}^{\it op} \ar[d]^{\Phi'} \ar[rr]^{\cal Y} &&  \cat{Set}^{\cat{C'}}\ar[dll]^<{{\it Lan}_{\cal Y}\Phi'} &&& \cat{C''}^{\it op} \ar[lll]^{\Psi'}\ar[dlllll]^{(\Psi\comp\Phi)'} \\
\cat{Set}^{\cat{C}} &&&&&} 
$$


Now we take ``profunctor glasses'' to look at Kelly's and Joyal's diagrams.

\begin{itemize}
\item A presheaf $X:?\cat{1}\rightarrow\cat{Set}$ can be viewed as a profunctor
$$(n\mapsto( \cdot\mapsto Xn)):?\cat{1}\profarrow\cat{1}$$

\item $X^{\otimes\_}$ can be viewed as a profunctor from $?\cat{1}$ to $?\cat{1}$ (since $(?\cat{1})^{\it op}=!(\cat{1}^{\it op})=!\cat{1}$). Thus the transformation
$$\seq{X:?\cat{1}\profarrow\cat{1}}{X^{\otimes\_}:?\cat{1}\profarrow?\cat{1}}$$
suggests to consider ? as a comonad over \cat{Prof}.
\item The diagram and the formula defining $Y\bullet X$ in Section \ref{Kelly-op} exhibit
$Y\bullet X$ as the profunctor composition of $Y:?\cat{1}\profarrow\cat{1}$ and
$X^{\otimes\_}:?\cat{1}\profarrow ?\cat{1}$.
But a better way to put it is that $Y\bullet X$ is the composition of $Y$ and $X$ in the coKleisli bicategory
$\cat{Prof}_?$, and we shall see that it is indeed the case.
\item  Joyal's construction amounts to ``taking points". We first observe that 
$?\cat{\emptyset}=\cat{1}$, where $\cat{\emptyset}$ is the empty category. (This holds for any of the eight choices for $!$.) It follows that $\cat{\emptyset}$ is terminal in $\cat{Prof}_?$. In this category, the points of an object $\cat{C}$ are
the morphisms from $\emptyset$ to \cat{C}, i.e., the profunctors from \cat{1} to \cat{C}, i.e. the presheaves over $\cat{C}$.  In particular, the points of $\cat{1}$ are just sets. 
We shall see that 
$$\cat{Prof}_?[\emptyset,X]={\it Lan}_\subseteq X$$
i.e., that the analytic functor associated with the presheaf $X$ describes its pointwise behaviour in
$\cat{Prof}_?$.
\end{itemize}

Our goal is thus to figure out ? as a comonad on \cat{Prof}. We shall do this in three steps:
\begin{itemize}
\item In Section \ref{ProKlei}, we show that profunctors arise as a  Kleisli category for the presheaf construction $\cat{C}\mapsto\psh{\cat{C}}$.
\item In Section \ref{Dist-op}, we show that all the monads ! of Section \ref{Burroni-eq} distribute over {\it Psh}. Distributive laws are recalled in Section \ref{Distrib-law}.
In Section \ref{intermezzo-sec}, we pause to explain how both Day's tensor product and this distributive law can be synthesised out of considerations of structure preservation.
\item in Section \ref{Putting-together}, we show that
this distributive law allows us to extend ! to a monad on \cat{Prof}, and by self-duality we obtain the comonad ? that we are looking for.
\end{itemize}

\smallskip
\paragraph{Warning.} In what follows, we ignore coherence issues for simplicity, and we partially address size issues.
\begin{itemize}
\item
Coherence issues arise from the fact that we shall compose morphisms using coend formulas, which make sense only up to iso. In particular, {\em our ``distributive law'' will in fact be a ``pseudo-distributive law''}~[\cite{Marmolejo,CHP,Gambino-Coh}]~.
\item 
Size issues arise from the fact that the presheaf construction is dramatically size increasing.
It is therefore in fact simply not rigorous to call {\it Psh} a monad (or even a pseudo-monad) on \cat{Cat}. 
But the issue is fortunately not too severe.
In a forthcoming paper,
Fiore, Gambino, Hyland, and Winskel propose a general notion of Kleisli structure,
 which we shall sketch here (ignoring again coherence issues), and in which {\it Psh} fits. 
 \end{itemize}

 \section{Profunctors as a Kleisli category} \label{ProKlei}
If we write a profunctor $\Phi:\cat{C_1}\profarrow \cat{C_2}$ in (untwisted) curried form
$$ C_1 \mapsto (C_2 \mapsto \Phi (C_1,C_2)):\cat{C_1}\rightarrow\cat{Set}^{\cat{C_2}^{\it op}}$$
this suggests us to look at the operation 
${\it Psh}$ on the objects of $\cat{Cat}$ defined by ${\it Psh}(\cat{C})=\cat{Set}^{\cat{C}^{\it op}}$.
The idea is to exhibit {\it Psh} as a (pseudo-)monad, 
so that profunctors arise as a Kleisli category. We have a good candidate for the
unit $\eta$, namely: ${\cal Y}:\cat{C}\rightarrow{\it Psh}(\cat{C})$.

But the Yoneda functor makes sense only for a locally small category (one in which all homsets are sets), while it is not clear at all whether {\it Psh} keeps us within the realm of locally small categories, i.e. if $\cat{C}$ is locally small, then $\psh{C}$ is not necessarily locally small.  
It is neverheless tempting to go on with  a multiplication $\mu:\psh{(\psh{C})}\rightarrow\psh{C}$ given by a left Kan extension
$$
 \xymatrix @C=0.4cm  {\psh{C} \ar[d]_{\it id} \ar[rr]^{\cal Y} &&  \psh{(\psh{C})}\ar[dll]^{\;\;\mu={\it Lan}_{\cal Y} ({\it id})}\\
\psh{C} &&}
$$
with explicit formula $\mu(H)C=\int^F HF\times FC$.
But why should a coend on such a vertiginous indexing collection exist?

\smallskip
We pause here to recall that an equivalent presentation of a monad on a category $\cat{C}$ is by means of  the following data:
\begin{itemize}
\item for each object $C$ of $\cat{C}$, an object $TC$ of \cat{C}, and a morphism $\eta_C\in\cat{C}[C,TC]$;
\item for all objects $C,D$ of \cat{C} and each morphism $f\in\cat{C}[C,TD]$, a morphism $f^{\#}\in\cat{C}[TC,TD]$;
\end{itemize} 
satisfying the equations $f=f^{\#}\comp \eta=f,\;,\eta^{\#}={\it id},\;,(g^{\#}\comp f)^{\#}=g^{\#}\comp f^{\#}$. 

Let us also recall the definition of the associated Kleisli category $\cat{C}_T$: its objects are 
the objects of \cat{C}, and one sets $\cat{C}_T[C,D]=\cat{C}[C,TD]$, with composition easily defined using the composition in \cat{C} and the $\_^{\#}$ operation.

\smallskip  It turns out that under this guise, the definition of monad can be generalised in a way that will fit our purposes.  Ignoring coherence issues, a {\em Kleisli structure} (as proposed in 
~[\cite{FGHW-Kleisli}]~) is given by the following data:
\begin{itemize}
\item a collection ${\cal A}$  of objects of \cat{C};
\item for each object $A$ in ${\cal A}$ an object $TA$ of \cat{C}, and a morphism $\eta_A\in\cat{C}[A,TA]$;
\item for all objects $A,B$ in ${\cal A}$ and each morphism $f\in\cat{C}[A,TB]$, a morphism $f^{\#}\in\cat{C}[TA,TB]$;
\end{itemize} 
satisfying the same equations.   The associated Kleisli category $\cat{C}_T$ has ${\cal A}$ as collection of objects, 
 and one sets $\cat{C}_T[A,B]=\cat{C}[A,TB]$. We recover the monads when  every object of \cat{C} is in ${\cal A}$.  But in general there is no such thing as a multiplication,  since it is not granted that $TTC$ exists for all $C$.  We also note that in a Kleisli structure, $T$ is still a functor, not from 
 \cat{C} to \cat{C}, but from $\cat{C}\downharpoonleft_{{\cal A}}$ (the full subcategory spanned by 
 ${\cal A}$) to \cat{C}.
 
 In our setting, we have $\cat{C}=\cat{Cat}$, and we can take ${\cal A}$ to consist of the small categories (in which the objects form a set, as well as all homsets). Then we have no worry about the coends that we shall write.  We complete the definition of the Kleisli structure with the definition of the $\_^{\#}$ operation.
 $$
 \xymatrix @C=0.4cm  {\cat{C} \ar[d]_{F} \ar[r]^{\cal Y} &  \cat{Set}^{\cat{C}^{\it op}}\ar[dl]^{\;\;\;F^{\#}={\it Lan}_{\cal Y} (F)}\\
\cat{Set}^{\cat{D}^{\it op}} &} 
$$
 and the explicit formula is  $F^{\#}XD=\int^C XC\times FCD$.
 It is then an easy exercise to check that the three equations of a Kleisli structure are satisfied, and that the composition in the associated Kleisli category coincides with the composition of profunctors as defined in the previous section.

We end the section with a description of the functorial action of ${\it Psh}$ (on, say, functors between small categories):
 $${\it Psh}(F)XD=(\eta_{\it Psh}\comp F)^{\#}XD=\int^C XC\times\cat{D}[D,FC]$$ 

\section{Distributive laws} \label{Distrib-law}

Recall that a distributive law~[\cite{Beck}]~(see also~[\cite{TTT}]~) is a natural transformation
$\dist:TS\rightarrow ST$, where $(S,\eta_S,\mu_S)$ and $(T,\eta_T,\mu_T)$ are two monads over the same category \cat{C}, satisfying the following laws, expressed in the language of string diagrams:

\medskip
\noindent EQUATION $(\dist-\eta_T)$:
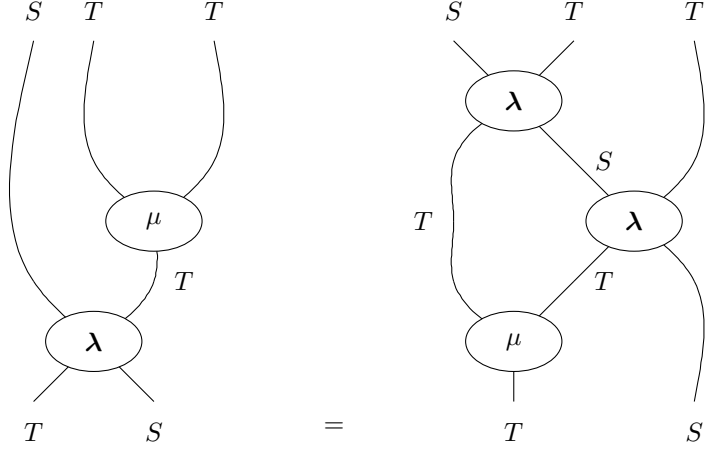
 $$\begin{array}{lllll}
\begin{tikzpicture}[scale=0.80]
\useasboundingbox (-0.5,-0.5) rectangle (2.5,3.5);
\draw (0.00,3.00) -- (-0.01,2.95) -- (-0.01,2.89) -- (-0.02,2.84) -- (-0.03,2.79) -- (-0.03,2.74) -- (-0.04,2.68) -- (-0.04,2.63) -- (-0.05,2.58) -- (-0.05,2.53) -- (-0.05,2.48) -- (-0.06,2.43) -- (-0.06,2.38) -- (-0.05,2.33) -- (-0.05,2.28) -- (-0.05,2.23) -- (-0.04,2.18) -- (-0.03,2.14) -- (-0.03,2.09) -- (-0.01,2.04) -- (0.00,2.00);
\draw (0.00,2.00) -- (0.02,1.94) -- (0.05,1.88) -- (0.08,1.82) -- (0.12,1.77) -- (0.15,1.71) -- (0.20,1.66) -- (0.24,1.60) -- (0.29,1.55) -- (0.34,1.50) -- (0.39,1.45) -- (0.44,1.41) -- (0.50,1.36) -- (0.56,1.31) -- (0.62,1.27) -- (0.68,1.22) -- (0.74,1.18) -- (0.81,1.13) -- (0.87,1.09) -- (0.94,1.04) -- (1.00,1.00);
\draw (2.00,3.00) -- (2.01,2.95) -- (2.01,2.89) -- (2.02,2.84) -- (2.03,2.79) -- (2.03,2.74) -- (2.04,2.68) -- (2.04,2.63) -- (2.05,2.58) -- (2.05,2.53) -- (2.05,2.48) -- (2.06,2.43) -- (2.06,2.38) -- (2.05,2.33) -- (2.05,2.28) -- (2.05,2.23) -- (2.04,2.18) -- (2.03,2.14) -- (2.03,2.09) -- (2.01,2.04) -- (2.00,2.00);
\draw (2.00,2.00) -- (1.98,1.94) -- (1.95,1.88) -- (1.92,1.82) -- (1.88,1.77) -- (1.85,1.71) -- (1.80,1.66) -- (1.76,1.60) -- (1.71,1.55) -- (1.66,1.50) -- (1.61,1.45) -- (1.56,1.41) -- (1.50,1.36) -- (1.44,1.31) -- (1.38,1.27) -- (1.32,1.22) -- (1.26,1.18) -- (1.19,1.13) -- (1.13,1.09) -- (1.06,1.04) -- (1.00,1.00);
\draw (1.00,1.00) -- (2.00,0.00);
\draw (1.00,1.00) -- (0.00,0.00);
\filldraw[fill=white] (2.00,3.00) ellipse (0.80cm and 0.50cm);
\filldraw[fill=white] (1.00,1.00) ellipse (0.80cm and 0.50cm);
\draw (0.00,3.50) node{$S$};
\draw (2.00,3.00) node{$\eta$};
\draw (2.50,2.00) node{$T$};
\draw (1.00,1.00) node{$\dist$};
\draw (0.00,-0.50) node{$T$};
\draw (2.00,-0.50) node{$S$};
\end{tikzpicture}
&\quad\quad& = &\quad\quad&  
\begin{tikzpicture}[scale=0.80]
\useasboundingbox (-0.5,-0.5) rectangle (1.5,3.5);
\draw (0.00,3.00) -- (0.00,0.00);
\draw (1.00,3.50) -- (1.00,3.47) -- (1.00,3.45) -- (1.00,3.42) -- (1.00,3.40) -- (1.00,3.38) -- (1.00,3.35) -- (1.00,3.33) -- (1.00,3.30) -- (1.00,3.27) -- (1.00,3.25) -- (1.00,3.23) -- (1.00,3.20) -- (1.00,3.17) -- (1.00,3.15) -- (1.00,3.12) -- (1.00,3.10) -- (1.00,3.07) -- (1.00,3.05) -- (1.00,3.02) -- (1.00,3.00);
\draw (1.00,3.00) -- (1.00,2.82) -- (1.00,2.65) -- (1.00,2.48) -- (1.00,2.30) -- (1.00,2.12) -- (1.00,1.95) -- (1.00,1.77) -- (1.00,1.60) -- (1.00,1.42) -- (1.00,1.25) -- (1.00,1.08) -- (1.00,0.90) -- (1.00,0.73) -- (1.00,0.55) -- (1.00,0.38) -- (1.00,0.20) -- (1.00,0.03) -- (1.00,-0.15) -- (1.00,-0.32) -- (1.00,-0.50);
\filldraw[fill=white] (0.00,3.00) ellipse (0.80cm and 0.50cm);
\draw (0.00,3.00) node{$\eta$};
\draw (1.50,2.00) node{$S$};
\draw (0.00,-0.50) node{$T$};
\end{tikzpicture}
\end{array}$$

\medskip
 \noindent EQUATION $(\dist-\mu_T)$:
$$\begin{array}{lllll}
\begin{tikzpicture}[scale=0.80]
\useasboundingbox (-0.5,-0.5) rectangle (3.5,6.5);
\draw (0.00,6.00) -- (-0.05,5.78) -- (-0.10,5.56) -- (-0.15,5.34) -- (-0.20,5.12) -- (-0.24,4.90) -- (-0.29,4.68) -- (-0.32,4.47) -- (-0.35,4.25) -- (-0.37,4.04) -- (-0.39,3.84) -- (-0.40,3.63) -- (-0.40,3.43) -- (-0.39,3.24) -- (-0.37,3.05) -- (-0.34,2.86) -- (-0.30,2.68) -- (-0.25,2.50) -- (-0.18,2.33) -- (-0.10,2.16) -- (0.00,2.00);
\draw (0.00,2.00) -- (0.04,1.94) -- (0.08,1.89) -- (0.12,1.84) -- (0.16,1.78) -- (0.21,1.73) -- (0.25,1.68) -- (0.30,1.63) -- (0.35,1.58) -- (0.40,1.53) -- (0.45,1.48) -- (0.50,1.43) -- (0.56,1.38) -- (0.61,1.33) -- (0.66,1.29) -- (0.72,1.24) -- (0.77,1.19) -- (0.83,1.14) -- (0.89,1.09) -- (0.94,1.05) -- (1.00,1.00);
\draw (1.00,6.00) -- (0.98,5.89) -- (0.96,5.78) -- (0.94,5.67) -- (0.92,5.57) -- (0.90,5.46) -- (0.89,5.35) -- (0.87,5.25) -- (0.86,5.14) -- (0.85,5.04) -- (0.84,4.94) -- (0.84,4.83) -- (0.84,4.73) -- (0.84,4.64) -- (0.85,4.54) -- (0.86,4.44) -- (0.88,4.35) -- (0.90,4.26) -- (0.93,4.17) -- (0.96,4.08) -- (1.00,4.00);
\draw (1.00,4.00) -- (1.03,3.94) -- (1.06,3.89) -- (1.10,3.83) -- (1.14,3.78) -- (1.18,3.72) -- (1.23,3.67) -- (1.27,3.62) -- (1.32,3.57) -- (1.37,3.52) -- (1.42,3.47) -- (1.48,3.42) -- (1.53,3.37) -- (1.59,3.32) -- (1.64,3.28) -- (1.70,3.23) -- (1.76,3.18) -- (1.82,3.14) -- (1.88,3.09) -- (1.94,3.05) -- (2.00,3.00);
\draw (3.00,6.00) -- (3.02,5.89) -- (3.04,5.78) -- (3.06,5.67) -- (3.08,5.57) -- (3.10,5.46) -- (3.11,5.35) -- (3.13,5.25) -- (3.14,5.14) -- (3.15,5.04) -- (3.16,4.94) -- (3.16,4.83) -- (3.16,4.73) -- (3.16,4.64) -- (3.15,4.54) -- (3.14,4.44) -- (3.12,4.35) -- (3.10,4.26) -- (3.07,4.17) -- (3.04,4.08) -- (3.00,4.00);
\draw (3.00,4.00) -- (2.97,3.94) -- (2.94,3.89) -- (2.90,3.83) -- (2.86,3.78) -- (2.82,3.72) -- (2.77,3.67) -- (2.73,3.62) -- (2.68,3.57) -- (2.63,3.52) -- (2.58,3.47) -- (2.52,3.42) -- (2.47,3.37) -- (2.41,3.32) -- (2.36,3.28) -- (2.30,3.23) -- (2.24,3.18) -- (2.18,3.14) -- (2.12,3.09) -- (2.06,3.05) -- (2.00,3.00);
\draw (2.00,3.00) -- (2.01,2.95) -- (2.01,2.89) -- (2.02,2.84) -- (2.03,2.79) -- (2.03,2.74) -- (2.04,2.68) -- (2.04,2.63) -- (2.05,2.58) -- (2.05,2.53) -- (2.05,2.48) -- (2.06,2.43) -- (2.06,2.38) -- (2.05,2.33) -- (2.05,2.28) -- (2.05,2.23) -- (2.04,2.18) -- (2.03,2.14) -- (2.03,2.09) -- (2.01,2.04) -- (2.00,2.00);
\draw (2.00,2.00) -- (1.98,1.94) -- (1.95,1.88) -- (1.92,1.82) -- (1.88,1.77) -- (1.85,1.71) -- (1.80,1.66) -- (1.76,1.60) -- (1.71,1.55) -- (1.66,1.50) -- (1.61,1.45) -- (1.56,1.41) -- (1.50,1.36) -- (1.44,1.31) -- (1.38,1.27) -- (1.32,1.22) -- (1.26,1.18) -- (1.19,1.13) -- (1.13,1.09) -- (1.06,1.04) -- (1.00,1.00);
\draw (1.00,1.00) -- (2.00,0.00);
\draw (1.00,1.00) -- (0.00,0.00);
\filldraw[fill=white] (2.00,3.00) ellipse (0.80cm and 0.50cm);
\filldraw[fill=white] (1.00,1.00) ellipse (0.80cm and 0.50cm);
\draw (0.00,6.50) node{$S$};
\draw (1.00,6.50) node{$T$};
\draw (3.00,6.50) node{$T$};
\draw (2.00,3.00) node{$\mu$};
\draw (2.50,2.00) node{$T$};
\draw (1.00,1.00) node{$\dist$};
\draw (0.00,-0.50) node{$T$};
\draw (2.00,-0.50) node{$S$};
\end{tikzpicture}
&\quad\quad& = &\quad\quad&  
\begin{tikzpicture}[scale=0.80]
\useasboundingbox (-0.5,-0.5) rectangle (5.5,6.5);
\draw (4.00,6.00) -- (4.02,5.89) -- (4.04,5.78) -- (4.06,5.67) -- (4.08,5.57) -- (4.10,5.46) -- (4.11,5.35) -- (4.13,5.25) -- (4.14,5.14) -- (4.15,5.04) -- (4.16,4.94) -- (4.16,4.83) -- (4.16,4.73) -- (4.16,4.64) -- (4.15,4.54) -- (4.14,4.44) -- (4.12,4.35) -- (4.10,4.26) -- (4.07,4.17) -- (4.04,4.08) -- (4.00,4.00);
\draw (4.00,4.00) -- (3.97,3.94) -- (3.94,3.89) -- (3.90,3.83) -- (3.86,3.78) -- (3.82,3.72) -- (3.77,3.67) -- (3.73,3.62) -- (3.68,3.57) -- (3.63,3.52) -- (3.58,3.47) -- (3.52,3.42) -- (3.47,3.37) -- (3.41,3.32) -- (3.36,3.28) -- (3.30,3.23) -- (3.24,3.18) -- (3.18,3.14) -- (3.12,3.09) -- (3.06,3.05) -- (3.00,3.00);
\draw (1.00,5.00) -- (1.05,4.95) -- (1.10,4.90) -- (1.15,4.85) -- (1.20,4.80) -- (1.25,4.75) -- (1.30,4.70) -- (1.35,4.65) -- (1.40,4.60) -- (1.45,4.55) -- (1.50,4.50) -- (1.55,4.45) -- (1.60,4.40) -- (1.65,4.35) -- (1.70,4.30) -- (1.75,4.25) -- (1.80,4.20) -- (1.85,4.15) -- (1.90,4.10) -- (1.95,4.05) -- (2.00,4.00);
\draw (2.00,4.00) -- (2.05,3.95) -- (2.10,3.90) -- (2.15,3.85) -- (2.20,3.80) -- (2.25,3.75) -- (2.30,3.70) -- (2.35,3.65) -- (2.40,3.60) -- (2.45,3.55) -- (2.50,3.50) -- (2.55,3.45) -- (2.60,3.40) -- (2.65,3.35) -- (2.70,3.30) -- (2.75,3.25) -- (2.80,3.20) -- (2.85,3.15) -- (2.90,3.10) -- (2.95,3.05) -- (3.00,3.00);
\draw (1.00,5.00) -- (0.93,4.96) -- (0.87,4.91) -- (0.80,4.87) -- (0.74,4.82) -- (0.67,4.78) -- (0.61,4.73) -- (0.55,4.69) -- (0.49,4.64) -- (0.43,4.59) -- (0.38,4.55) -- (0.32,4.50) -- (0.27,4.45) -- (0.23,4.40) -- (0.18,4.34) -- (0.14,4.29) -- (0.11,4.23) -- (0.07,4.18) -- (0.04,4.12) -- (0.02,4.06) -- (0.00,4.00);
\draw (0.00,4.00) -- (-0.01,3.96) -- (-0.02,3.91) -- (-0.03,3.86) -- (-0.03,3.82) -- (-0.03,3.77) -- (-0.04,3.72) -- (-0.04,3.67) -- (-0.04,3.62) -- (-0.03,3.57) -- (-0.03,3.52) -- (-0.03,3.47) -- (-0.02,3.42) -- (-0.02,3.37) -- (-0.02,3.32) -- (-0.01,3.26) -- (-0.01,3.21) -- (-0.00,3.16) -- (-0.00,3.11) -- (-0.00,3.05) -- (0.00,3.00);
\draw (0.00,3.00) -- (-0.00,2.95) -- (-0.00,2.89) -- (-0.00,2.84) -- (-0.01,2.79) -- (-0.01,2.74) -- (-0.02,2.68) -- (-0.02,2.63) -- (-0.02,2.58) -- (-0.03,2.53) -- (-0.03,2.48) -- (-0.03,2.43) -- (-0.04,2.38) -- (-0.04,2.33) -- (-0.04,2.28) -- (-0.03,2.23) -- (-0.03,2.18) -- (-0.03,2.14) -- (-0.02,2.09) -- (-0.01,2.04) -- (0.00,2.00);
\draw (0.00,2.00) -- (0.02,1.94) -- (0.04,1.88) -- (0.07,1.82) -- (0.11,1.77) -- (0.14,1.71) -- (0.18,1.66) -- (0.23,1.60) -- (0.27,1.55) -- (0.32,1.50) -- (0.38,1.45) -- (0.43,1.41) -- (0.49,1.36) -- (0.55,1.31) -- (0.61,1.27) -- (0.67,1.22) -- (0.74,1.18) -- (0.80,1.13) -- (0.87,1.09) -- (0.93,1.04) -- (1.00,1.00);
\draw (2.00,6.00) -- (1.00,5.00);
\draw (0.00,6.00) -- (1.00,5.00);
\draw (3.00,3.00) -- (3.06,2.95) -- (3.12,2.91) -- (3.18,2.86) -- (3.24,2.82) -- (3.30,2.77) -- (3.36,2.72) -- (3.41,2.68) -- (3.47,2.63) -- (3.52,2.58) -- (3.58,2.53) -- (3.63,2.48) -- (3.68,2.43) -- (3.73,2.38) -- (3.77,2.33) -- (3.82,2.28) -- (3.86,2.22) -- (3.90,2.17) -- (3.94,2.11) -- (3.97,2.06) -- (4.00,2.00);
\draw (4.00,2.00) -- (4.04,1.92) -- (4.07,1.83) -- (4.10,1.74) -- (4.12,1.65) -- (4.14,1.56) -- (4.15,1.46) -- (4.16,1.36) -- (4.16,1.27) -- (4.16,1.17) -- (4.16,1.06) -- (4.15,0.96) -- (4.14,0.86) -- (4.13,0.75) -- (4.11,0.65) -- (4.10,0.54) -- (4.08,0.43) -- (4.06,0.33) -- (4.04,0.22) -- (4.02,0.11) -- (4.00,0.00);
\draw (3.00,3.00) -- (2.95,2.95) -- (2.90,2.90) -- (2.85,2.85) -- (2.80,2.80) -- (2.75,2.75) -- (2.70,2.70) -- (2.65,2.65) -- (2.60,2.60) -- (2.55,2.55) -- (2.50,2.50) -- (2.45,2.45) -- (2.40,2.40) -- (2.35,2.35) -- (2.30,2.30) -- (2.25,2.25) -- (2.20,2.20) -- (2.15,2.15) -- (2.10,2.10) -- (2.05,2.05) -- (2.00,2.00);
\draw (2.00,2.00) -- (1.95,1.95) -- (1.90,1.90) -- (1.85,1.85) -- (1.80,1.80) -- (1.75,1.75) -- (1.70,1.70) -- (1.65,1.65) -- (1.60,1.60) -- (1.55,1.55) -- (1.50,1.50) -- (1.45,1.45) -- (1.40,1.40) -- (1.35,1.35) -- (1.30,1.30) -- (1.25,1.25) -- (1.20,1.20) -- (1.15,1.15) -- (1.10,1.10) -- (1.05,1.05) -- (1.00,1.00);
\draw (1.00,1.00) -- (1.00,0.00);
\filldraw[fill=white] (1.00,5.00) ellipse (0.80cm and 0.50cm);
\filldraw[fill=white] (3.00,3.00) ellipse (0.80cm and 0.50cm);
\filldraw[fill=white] (1.00,1.00) ellipse (0.80cm and 0.50cm);
\draw (0.00,6.50) node{$S$};
\draw (2.00,6.50) node{$T$};
\draw (4.00,6.50) node{$T$};
\draw (1.00,5.00) node{$\dist$};
\draw (2.50,4.00) node{$S$};
\draw (-0.50,3.00) node{$T$};
\draw (3.00,3.00) node{$\dist$};
\draw (2.50,2.00) node{$T$};
\draw (1.00,1.00) node{$\mu$};
\draw (1.00,-0.50) node{$T$};
\draw (4.00,-0.50) node{$S$};
\end{tikzpicture}
\end{array}$$
and two similar equations $(\dist-\mu_S)$ and $(\dist-\eta_S)$.

A distributive law allows us to extend $T$ to $\cat{C}_S$, the Kleisli category of $S$ (and conversely, such a lifting induces a distributive law such that the two constructions are inverse to each other).\footnote{\label{dist-lift-alg} From another standpoint, a distributive law $\dist:TS\rightarrow ST$ also amounts to lifting $S$ to the category of $T$-algebras. See Section \ref{intermezzo-sec}.}
The extended $T$ acts on objects as the old $T$. On morphisms, its action is described as follows:
$$
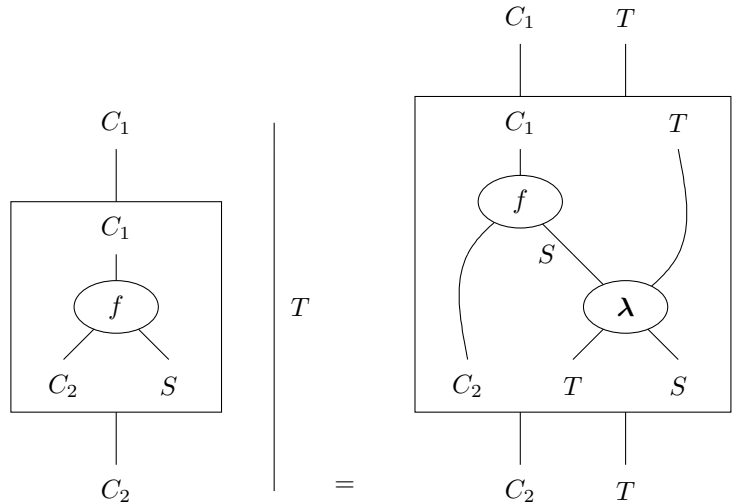
\begin{tikzpicture}[scale=0.7]
\useasboundingbox (-0.5,-0.5) rectangle (5.5,6.5);
\draw (5.00,6.50) -- (5.00,6.47) -- (5.00,6.45) -- (5.00,6.42) -- (5.00,6.40) -- (5.00,6.38) -- (5.00,6.35) -- (5.00,6.33) -- (5.00,6.30) -- (5.00,6.28) -- (5.00,6.25) -- (5.00,6.22) -- (5.00,6.20) -- (5.00,6.17) -- (5.00,6.15) -- (5.00,6.12) -- (5.00,6.10) -- (5.00,6.08) -- (5.00,6.05) -- (5.00,6.02) -- (5.00,6.00);
\draw (5.00,6.00) -- (5.00,5.67) -- (5.00,5.35) -- (5.00,5.02) -- (5.00,4.70) -- (5.00,4.38) -- (5.00,4.05) -- (5.00,3.72) -- (5.00,3.40) -- (5.00,3.07) -- (5.00,2.75) -- (5.00,2.42) -- (5.00,2.10) -- (5.00,1.77) -- (5.00,1.45) -- (5.00,1.12) -- (5.00,0.80) -- (5.00,0.47) -- (5.00,0.15) -- (5.00,-0.18) -- (5.00,-0.50);
\draw (2.00,6.00) -- (2.00,5.00);
\draw (2.00,3.00) -- (3.00,2.00);
\draw (2.00,3.00) -- (1.00,2.00);
\draw (2.00,4.00) -- (2.00,3.00);
\draw (2.00,1.00) -- (2.00,0.00);
\draw[solid,] (0.00,5.00) rectangle (4.00,1.00);
\filldraw[fill=white] (2.00,3.00) ellipse (0.80cm and 0.50cm);
\draw (2.00,6.50) node{$C_1$};
\draw (2.00,4.50) node{$C_1$};
\draw (2.00,3.00) node{$f$};
\draw (5.50,3.00) node{$T$};
\draw (1.00,1.50) node{$C_2$};
\draw (3.00,1.50) node{$S$};
\draw (2.00,-0.50) node{$C_2$};
\end{tikzpicture}
\quad = \quad 
\begin{tikzpicture}[scale=0.7]
\useasboundingbox (-0.5,-0.5) rectangle (6.5,8.5);
\draw (2.00,8.00) -- (2.00,7.00);
\draw (4.00,8.00) -- (4.00,7.00);
\draw (5.00,6.00) -- (5.02,5.89) -- (5.04,5.78) -- (5.06,5.67) -- (5.08,5.57) -- (5.10,5.46) -- (5.11,5.35) -- (5.13,5.25) -- (5.14,5.14) -- (5.15,5.04) -- (5.16,4.94) -- (5.16,4.83) -- (5.16,4.73) -- (5.16,4.64) -- (5.15,4.54) -- (5.14,4.44) -- (5.12,4.35) -- (5.10,4.26) -- (5.07,4.17) -- (5.04,4.08) -- (5.00,4.00);
\draw (5.00,4.00) -- (4.97,3.94) -- (4.94,3.89) -- (4.90,3.83) -- (4.86,3.78) -- (4.82,3.72) -- (4.77,3.67) -- (4.73,3.62) -- (4.68,3.57) -- (4.63,3.52) -- (4.58,3.47) -- (4.52,3.42) -- (4.47,3.37) -- (4.41,3.32) -- (4.36,3.28) -- (4.30,3.23) -- (4.24,3.18) -- (4.18,3.14) -- (4.12,3.09) -- (4.06,3.05) -- (4.00,3.00);
\draw (2.00,5.00) -- (2.05,4.95) -- (2.10,4.90) -- (2.15,4.85) -- (2.20,4.80) -- (2.25,4.75) -- (2.30,4.70) -- (2.35,4.65) -- (2.40,4.60) -- (2.45,4.55) -- (2.50,4.50) -- (2.55,4.45) -- (2.60,4.40) -- (2.65,4.35) -- (2.70,4.30) -- (2.75,4.25) -- (2.80,4.20) -- (2.85,4.15) -- (2.90,4.10) -- (2.95,4.05) -- (3.00,4.00);
\draw (3.00,4.00) -- (3.05,3.95) -- (3.10,3.90) -- (3.15,3.85) -- (3.20,3.80) -- (3.25,3.75) -- (3.30,3.70) -- (3.35,3.65) -- (3.40,3.60) -- (3.45,3.55) -- (3.50,3.50) -- (3.55,3.45) -- (3.60,3.40) -- (3.65,3.35) -- (3.70,3.30) -- (3.75,3.25) -- (3.80,3.20) -- (3.85,3.15) -- (3.90,3.10) -- (3.95,3.05) -- (4.00,3.00);
\draw (2.00,5.00) -- (1.94,4.95) -- (1.88,4.91) -- (1.82,4.86) -- (1.76,4.82) -- (1.70,4.77) -- (1.64,4.72) -- (1.59,4.68) -- (1.53,4.63) -- (1.48,4.58) -- (1.42,4.53) -- (1.37,4.48) -- (1.32,4.43) -- (1.27,4.38) -- (1.23,4.33) -- (1.18,4.28) -- (1.14,4.22) -- (1.10,4.17) -- (1.06,4.11) -- (1.03,4.06) -- (1.00,4.00);
\draw (1.00,4.00) -- (0.96,3.92) -- (0.93,3.83) -- (0.90,3.74) -- (0.88,3.65) -- (0.86,3.56) -- (0.85,3.46) -- (0.84,3.36) -- (0.84,3.27) -- (0.84,3.17) -- (0.84,3.06) -- (0.85,2.96) -- (0.86,2.86) -- (0.87,2.75) -- (0.89,2.65) -- (0.90,2.54) -- (0.92,2.43) -- (0.94,2.33) -- (0.96,2.22) -- (0.98,2.11) -- (1.00,2.00);
\draw (2.00,6.00) -- (2.00,5.00);
\draw (4.00,3.00) -- (5.00,2.00);
\draw (4.00,3.00) -- (3.00,2.00);
\draw (2.00,1.00) -- (2.00,0.00);
\draw (4.00,1.00) -- (4.00,0.00);
\draw[solid,] (0.00,7.00) rectangle (6.00,1.00);
\filldraw[fill=white] (2.00,5.00) ellipse (0.80cm and 0.50cm);
\filldraw[fill=white] (4.00,3.00) ellipse (0.80cm and 0.50cm);
\draw (2.00,8.50) node{$C_1$};
\draw (4.00,8.50) node{$T$};
\draw (2.00,6.50) node{$C_1$};
\draw (5.00,6.50) node{$T$};
\draw (2.00,5.00) node{$f$};
\draw (2.50,4.00) node{$S$};
\draw (4.00,3.00) node{$\dist$};
\draw (1.00,1.50) node{$C_2$};
\draw (3.00,1.50) node{$T$};
\draw (5.00,1.50) node{$S$};
\draw (2.00,-0.50) node{$C_2$};
\draw (4.00,-0.50) node{$T$};
\end{tikzpicture}
$$
where the box separates the inside which lives in \cat{C} from the outside which lives in $\cat{C}_S$.
It is a nice exercise to check that the old and the new $T$ satisfy
$$T\comp F_S=F_S\comp T$$
where $F_S:\cat{C}\rightarrow\cat{C}_S$ is defined by $F_S(C)=C$ and $F_S(f)=\eta_S\comp f$,
and that the new $T$ is a monad.

\smallskip
In the next section, we want to apply this to $S={\it Psh}$ and $T=!$ .  But $S$ is only a Kleisli structure.  Only slight adjustments is needed: 
\begin{enumerate}
\item make sure that  ${\cal A}$ is stable under $T$; 
\item replace the equation  $(\dist-\mu_S)$ by the following one  (for $f:A\rightarrow SB$):
$$\dist_B\comp T(g^{\#})=(\dist_B\comp Tg)^{\#}\comp\dist_A$$
\end{enumerate}

\section{A $!/{\it Psh}$ distributive law} \label{Dist-op}
We define a transformation $\dist:\:!\comp{\it Psh}\rightarrow{\it Psh}\:\comp \:!$, as follows:

\Proofitem{\bullet} on objects:
\begin{center}
\fbox{$\begin{array}{l}(F_1,\ldots,F_n)\mapsto ((C_1,\ldots,C_m) \mapsto \\
 \quad\quad \quad\quad 
 \int^{D_1,\ldots,D_n} F_1D_1\times\ldots\times
 F_n D_n\times !\cat{C}[(C_1,\ldots,C_m),(D_1,\ldots,D_n)])
 \end{array}$}
 \end{center}

\Proofitem{\bullet} on generating morphisms:
 \begin{itemize}
 \item[-] $\alpha:F\rightarrow G$. There is an obvious map from
$\int^D FD\times !\cat{C}[C_1,\ldots,C_m,(D)]$ to $\int^D GD\times !\cat{C}[C_1,\ldots,C_m,(D)]$.
\item[-] Combinators, say $\delta:(F)\rightarrow (F,F)$. We obtain a map
$$\begin{array}{l}
\mbox{from }
\int^D FD\times \mathord{!}\cat{C}[C_1,\ldots,C_m,(D)]\\
\mbox{to } 
\int^{D_1,D_2} FD_1\times
 F D_2\times \mathord{!}\cat{C}[(C_1,\ldots,C_m,(D_1,D_2)]
 \end{array}$$ 
 by going 
$$ \begin{array}{l}
 \mbox{from }FD\times \mathord{!}\cat{C}[C_1,\ldots,C_m,(D)]\\
 \mbox{to }FD\times FD\times \mathord{!}\cat{C}[(C_1,\ldots,C_m,(D,D)]
 \end{array}$$
\end{itemize}
We verify three of the four equations (leaving the last one -- as stated at the end of Section \ref{Distrib-law} -- to the reader):

\smallskip\noindent
 $(\dist-\eta_{!})$.  The left hand side is computed by taking the case $n=1$ in the definition of $\dist$, thus we obtain
$$F\mapsto ((C_1,\ldots,C_m)\mapsto \int^D FD\times !\cat{C}[(C_1,\ldots,C_m),(D)])$$
which is the formula for  ${\it Psh}(\eta_!)$.

\smallskip\noindent
 $(\dist-\eta_{\it Psh})$. We have (at $(A_1,\ldots,A_n)$, $(C_1,\ldots C_m)$):
$$\begin{array}{l}
\int^{D_1,\ldots,D_n}\cat{C}[D_1,A_1]\times\ldots\times
\cat{C}[D_n,A_n]\times!\cat{C}[(C_1,\ldots,C_m),(D_1,\ldots,D_n)]\\
\quad =  !\cat{C}[(C_1,\ldots,C_m),(A_1,\ldots,A_n)] ={\cal Y}(A_1,\ldots,A_n)(C_1,\ldots C_m)
\end{array}$$

\smallskip\noindent
 $(\dist-\mu_!)$. The left hand side is (at $((F_1^1,\ldots,F_1^{i_1}),\ldots(F_n^1,\ldots,F_n^{i_n}))$,
 $(C_1,\ldots,C_m)$):
 $$\int^{D_1^1,\ldots,D_n^{i_n}}F_1^1D_1^1\times\ldots\times F_n^{i_n}D_n^{i_n}\times!\cat{C}[\vec{C},(D_1^1,\ldots,D_n^{i_n})]$$
 We first compute the upper part $(\dist !)\comp (!\dist)$ of the right hand side (at $(\vec{F_1},\ldots,\vec{F_n})$, $((A_1^1,\ldots,A_1^{j_1}),\ldots,(A_p^1,\ldots,A_p^{j_p}))$):
 $$
 \int^{\vec{B_1},\ldots\vec{B_n}}\!\!\!\!\!\int^{D_1^1,\ldots,D_n^{i_n}}\!\!\!\!\!\!\!\!F_1^1D_1^1\times\cdots\times F_n^{i_n}D_n^{i_n}\times!\cat{C}[\vec{B_1},\vec{D_1}]\times\cdots\times
 !\cat{C}[\vec{C_n},\vec{D_n}]\times !!\cat{C}[\vec{\vec{A}},\vec{\vec{B}}]$$
$$= \int^{D_1^1,\ldots,D_n^{i_n}}F_1^1D_1^1\times\ldots\times F_n^{i_n}D_n^{i_n}\times !!\cat{C}[\vec{\vec{A}},\vec{\vec{D}}]\quad\quad\quad\quad\quad\quad
$$
 Finally, applying ${\it Psh}(\mu)$, we get (at $\vec{\vec{F}}$, $\vec{C}$):
 $$\int^{\vec{A_1},\ldots,\vec{A_p}}\!\!\!\int^{D_1^1,\ldots,D_n^{i_n}}\!\!\!\!F_1^1D_1^1\times\cdots\times F_n^{i_n}D_n^{i_n}\times
!!\cat{C}[\vec{\vec{A}},\vec{\vec{D}}] \times !\cat{C}[\vec{C},(A_1^1,\ldots,A_p^{j_p})]$$
and we conclude using $\mu_!$.


\medskip
Summarizing, we have the following result:

\begin{proposition}
The transformation $\dist$ defines a (pseudo)-distributive law for any choice among  the six $!_X$ monads over the presheaf Kleisli structure.
\end{proposition}

A similar proposition is proved in~[\cite{Tanaka-Power}]~: there, the existence of $\dist$ is derived from the fact $!$ and ${\it Psh}$ are both free constructions (the latter being the free cocompletion, i.e. the free category with all colimits) and commute in the sense that {\it Psh} lifts to
the category of $!$-algebras (cf. Section \ref{Burroni-eq}). For completeness, we sketch this more conceptual approach in the following section.

On the other hand, the advantage of a ``symbol-pushing'' proof like the one presented above is to highlight plainly  the  ``uniformity" in the choice of any combination $X$ of Burroni's combinators.

\section{Intermezzo} \label{intermezzo-sec}

The reader may have noticed that the formula defining  $\dist$ proposed in the last section ``looks like" the formula for Day's tensor product (cf. Section \ref{Kelly-op}). In this section, we actually synthesise  the latter formula from considerations of cocontinuity, and then the former one from considerations of monad lifting.  Since this section has merely an explanatory purpose, we limit ourselves here to the case $!=!_\emptyset$, and we disregard size issues.

 Cocontinuity is at the heart of the presheaf construction. It is well-known that $\psh{C}$ is the free cocomplete category over $\cat{C}$, and that the unique cocontinuous (i.e., colimit preserving) extension of a functor $F:\cat{C}\rightarrow\cat{D}$ (where $\cat{D}$ is cocomplete, i.e., has all colimits) is its left Kan extension
 ${\it Lan}_{\cal Y}(F)$:
 
$$
  \xymatrix @C=0.4cm  {\cat{C} \ar[d]^{F} \ar[rr]^{\cal Y} &&  \psh{C}\ar[dll]^{{\it Lan}_{{\cal Y}}F}\\
\cat{D} &} 
$$
In particular, the fact that  ${\it Lan}_{\cal Y}(F)$ preserves colimits follows from its being left adjoint to $(C\mapsto(D\mapsto \cat{D}[FC,D]))$ (easy check).

We could have defined the ``monad'' ${\it Psh}$ more economically and more conceptually as the monad arising from the adjunction between $\cat{Cat}$ and the category of cocomplete categories. The ``brute force" construction of last section is in any case a good exercise in coend computations.

Although we shall not need it here, it is worth pointing out that the adjunction is monadic, i.e., that the (pseudo-) ${\it Psh}$-algebras are actually the cocomplete categories (and hence that there can be at most one ${\it Psh}$-algebra structure on a given category). This can be proved using Beck's characterization of monadic adjunctions (Theorem 4.4.4 of~[\cite{Borceux2}]~).

Consider a slightly more general version of Day's tensor product than that given in Section \ref{Kelly-op},
with now some monoidal category $\cat{C}$ in place of $!\cat{1}$
$$(X\otimes Y)C=\int^{C_1,C_2} XC_1\times XC_2\times\cat{C}[C,C_1\otimes C_2]$$
(thus Day's tensor depends on the tensor on $\cat{C}$).
We show that this definition is entirely determined from the requirements that
\begin{itemize}
\item $\otimes$ extends the tensor product of $\cat{C}$ (via the Yoneda embedding), and
\item $\otimes$ is cocontinuous in each argument.
\end{itemize}
And, of course, Day's product does satisfy these requirements (easy proof left to the reader).
Recall that ${\cal Y}$ is dense, i.e., that every preseheaf is a colimit of representable presheaves (i.e. presheaves of the form $YC$), and, more precisely, that ${\it Lan}_{\cal Y}{\cal Y}={\it id}$:
$$X=\int^D XD\cdot{\cal Y}D = {\it colim}_{(D,z\in XD)}{\cal Y}D$$
Then, we must have
$$\begin{array}{llll}
X\otimes Y & = & ({\it colim}_{(D_1,z\in XD_1)}{\cal Y}D_1)\otimes Y\\
& = & {\it colim}_{(D_1,z\in XD_1)}({\cal Y}D_1\otimes Y) & (\mbox{by requirement 2})\\
& = & {\it colim}_{(D_1,z\in XD_1),(D_2,z\in XD_2)} ({\cal Y}D_1\otimes {\cal Y}D_2)  & (\mbox{by requirement 2})\\
& = &  {\it colim}_{(D_1,z\in XD_1),(D_2,z\in XD_2)} {\cal Y}(D_1\otimes D_2) & (\mbox{by requirement 1})\\
& = & \int^{D_1,D_2} (XD_1\times XD_2)\cdot {\cal Y}(D_1\otimes D_2)
\end{array}$$
In other words, Day's tensor is determined, as announced.

It is shown in~[\cite{Im-Kelly}]~that the (2-) adjunction between categories and cocomplete categories specialises to an adjunction between monoidal categories and cocomplete monoidal categories, and that $\psh{C}$, equipped with Day's tensor product, is the free such one. We content ourselves here with 
the key verification, namely that when $F:\cat{C}\rightarrow \cat{D}$ is monoidal, then so is its unique cocontinuous extension $\tilde{F}={\it Lan}_{\cal Y}F$.
We first compute $\tilde{F}X_1\otimes \tilde{F}X_2$ and $\tilde{F}(X_1\otimes X_2)$:
$$\begin{array}{llll}
\tilde{F}X_1\otimes \tilde{F}X_2 & = & \int^{C_1,C_2}(X_1C_1\times X_2C_2)\cdot(FC_1\otimes FC_2) \quad (\mbox{by cocontinuity})\\\\
\tilde{F}(X_1\otimes X_2) & = & \int^C \int^{C_1,C_2} X_1C_1\times X_2C_2\times\cat{C}[C,C_1\otimes C_2]\cdot FC\\
& = &  \int^{C_1,C_2} X_1C_1\times X_2C_2\cdot F(C_1\otimes C_2)
\end{array}$$
From this, it follows that the maps $FC_1\otimes FC_2\rightarrow F(C_1\otimes C_2)$ induce a map
$\tilde{F}X_1\otimes \tilde{F}X_2\rightarrow \tilde{F}(X_1\otimes X_2)$.

As a consequence ${\it Psh}$ lifts to  (pseudo-) $!$-algebras, i.e.
to monoidal categories. Given two monads $T,S$, a lifting of $S$ to the category of $T$-algebras consists of a monad $S'$ on the category of $T$-algebras such that $S,S'$ commute with the forgetful functor, and such that the unit and the multiplication of $S'$ are mapped by the forgetful functor  to the unit and the multiplication of $S$.
Here, we take $S={\it Psh}$ and $T=!$. The lifting $S'$ maps a monoidal structure $\otimes$ on \cat{C} to the corresponding Day's monoidal structure on ${\it Psh}(\cat{C})$.

We turn back to the general situation.  A lifting of $S$ to $T$-algebras is equivalent to giving a distributive law $\dist:TS\rightarrow ST$ (cf. Footnote \ref{dist-lift-alg}). One defines $\dist$ from the data of the lifting as follows:
$$\dist_C= S'(\mu_C)\comp TS\eta_C\quad\quad(\mu,\eta\mbox{ relative to }T)$$
(recall that $\mu_C$ is a $T$-algebra --  the free one).
This allows us to derive the distributive law from Day's tensor product.
 First, for $S={\it Psh}$ and $T=!$,
$S'(\mu_{\cat{C}})$ is the
Day's tensor product $(G_1,\ldots,G_n)\mapsto G_1\otimes\ldots\otimes G_n$ 
 associated with the (flattening) tensor product on $!\cat{C}$:
$$(G_1\otimes\ldots\otimes G_n)\vec{C}=
\int^{\vec{A_1},\ldots,\vec{A_n}}G_1\vec{A_1}\times\ldots\times G_n\vec{A_n}\times !\cat{C}[\vec{C},(A_1^1,\ldots,A_n^{i_n})]$$
Then $\dist$  is obtained by replacing $G_i$ by ${\it Psh}(\eta)(F_i)$:
$$\int^{\vec{A_1},\ldots,\vec{A_n}} z_1
\times\ldots\times z_n\times !\cat{C}[\vec{C},(A_1^1,\ldots,A_n^{i_n})]$$
with $z_1=\int^{D_1}!F_1D_1\times!\cat{C}[\vec{A_1},(D_1)]$, \ldots,
$z_n= \int^{D_n}F_nD_n\times!\cat{C}[\vec{A_n},(D_n)]$. 
This
simplifies to the formula given in Section \ref{Dist-op}.

 \section{The (bi)category $\cat{Prof}_?$}  \label{Putting-together}

We next consider the coKleisli bicategory $\cat{Prof}_?=(\cat{Cat}_{\it Psh})_?$. From the previous sections we know
that the monad ! on \cat{Cat} extends to a (pseudo-) monad ! on \cat{Prof}, which under the self-duality
of \cat{Prof} gives a (pseudo-) comonad ? on \cat{Prof}. Just as a monad gives rise to a Kleisli category, a comonad gives rise to a coKleisli category. Here, $\cat{Prof}_?$ has categories as objects, while its (1-)morphisms are defined by
$$\cat{Prof}_?[\cat{C},\cat{C'}]=\cat{Prof}[?\cat{C},\cat{C'}]$$
For $X=\set{\sigma}$, this is the category  generalised species of structures of Fiore, Gambino, Hyland, and Winskel~[\cite{FGHW}]~.  Also with  $X=\set{\sigma}$, the endomorphisms 
in $\cat{Prof}[?\cat{C},\cat{C}]$ are the $C$-profiles of Baez and Dolan~[\cite{BD98}]~.
\footnote{\label{BD-not-dist} Baez and Dolan use a variation to synthesise the composition of $\cat{C}$-signatures:  they note that
$\cat{Set}^{?\cat{C}}$ is a  free symmetric monoidal cocomplete construction (where the tensor preserves the colimits in each argument, cf. Section \ref{intermezzo-sec}), which allows them to identify
$C$-profiles with the endofunctors on $\cat{Set}^{?\cat{C}}$ that preserve tensor and colimits, and hence to inherit composition from usual composition of functors.} 
\medskip 
We are now in a position to round the circle and to show (cf. Theorem 4.1 of~[\cite{FGHW}]~for the case $X=\set{\sigma}$)
that the composition in this category  coincides with the substitution operation as defined in Section \ref{Kelly-op}. We decompose all the steps of the construction.

\medskip\noindent
(1)
Let $\Phi:\cat{C}\profarrow\cat{C'}$. Then $!\Phi:\:!\cat{C}\profarrow \:!\cat{C'}$ is (using $\dist$):
$$!\Phi((C_1,\ldots,C_n),\vec{C'})=\int^{D'_1,\ldots,D'_n}\Phi C_1D'_1\times\ldots\times \Phi C_nD'_n\times
!\cat{C'}[\vec{C'},\vec{D'}]$$

\medskip\noindent
$(1)^{\it op}$ We can then define $?\Phi=(!\Phi^{\it op})^{\it op}$:
$$?\Phi(\vec{C},(C'_1,\ldots,C'_p)=\int^{D_1,\ldots,D_p}\Phi D_1C'_1\times\ldots\times\Phi D_pC'_p\times
?\cat{C}[\vec{D},\vec{C}]$$

\medskip\noindent
(2) The multiplication of $!$  viewed as a monad on \cat{Prof} at \cat{C} is 
$$(((C_1^1,\ldots,C_1^{i_1}),\ldots,(C_n^1,\ldots,C_n^{i_n})),\vec{D})\mapsto\: !\cat{C}[\vec{D},(C_1^1,\ldots,C_n^{i_n})]:\;!!\cat{C}\profarrow!\cat{C})$$

\medskip\noindent
$(2)^{\it op}$ Hence the comultiplication $\delta:\;?\cat{C}\profarrow??\cat{C}$ of $?=(\cat{C}\mapsto (!(\cat{C}^{\it op})^{\it op})$ is
$$ ((\vec{C},((D_1^1,\ldots,D_1^{i_1}),\ldots,(D_n^1,\ldots,D_n^{i_n})))\mapsto
?\cat{C}[(D_1^1,\ldots,D_n^{i_n}),\vec{C}])$$

\medskip\noindent
(3) Let $\Phi:?\cat{C}\profarrow\cat{C'}$. Then
$\Phi^{\flat}= (?\Phi)\comp\delta:?\cat{C}\profarrow?\cat{C'}$ is given by the following explicit formula (on $\vec{C},(C'_1,\ldots,C'_p)$:
$$\int^{\vec{A_1},\ldots,\vec{A_m}}\!\!\!\!\!\!\!?\cat{C}[(A_1^1,\ldots,A_m^{j_m}),\vec{C}]\times(\int^{\vec{D_1},\ldots,\vec{D_p}}\!\!\!\!\!\Phi\vec{D_1}C'_1\times
\dots\times\Phi\vec{D_p}C'_p\times??\cat{C}[\vec{\vec{D}},\vec{\vec{A}}])$$
$$= \int^{\vec{D_1},\ldots,\vec{D_p}}\Phi\vec{D_1}C'_1\times
\dots\times\Phi\vec{D_p}C'_p\times?\cat{C}[(D_1^1,\ldots,D_p^{i_p}),\vec{C}]
\quad\quad\quad\quad\quad\quad\quad$$
The simplification comes from the functoriality of the multiplication of the monad $!$ on \cat{Cat} (we have 
$\mu_{\cat{C}^{\it op}}:!!\cat{C}^{\it op}\rightarrow!\cat{C}^{\it op}$, which we can view as a functor
from $??\cat{C}=(!!\cat{C}^{\it op})^{\it op}$ to $?\cat{C}=(!\cat{C}^{\it op})^{\it op}$).

\medskip\noindent
Finally, we can compose $\Phi:?\cat{C}\profarrow\cat{C'}$ and $\Psi:?\cat{C'}\profarrow\cat{C''}$:
$\Psi\bullet\Phi =\Psi \comp \Phi^{\flat}$.  Explicitly, $(\Psi\bullet\Phi)(\vec{C},C'')$ is
given by the following formula:
$$\int^{C'_1,\ldots,C'_p}\!\!\!\!\!\Psi(\vec{C'},C'')\times\int^{\vec{D_1},\ldots,\vec{D_p}}\!\!\!\!\!\Phi\vec{D_1}C'_1\times
\dots\times\Phi\vec{D_p}C'_p\times?\cat{C}[(D_1^1,\ldots,D_p^{i_p}),\vec{C}]$$

\medskip
As promised, we recover Kelly's (Section \ref{Kelly-op}) and Joyal's (Section \ref{Joyal-op}) settings as instances:
\begin{enumerate}
\item
When $\cat{C}=\cat{C'}=\cat{C''}=\cat{1}$, then $\Psi$ and $\Phi$ are presheaves $Y,X$, and the definition boils down to $Y\bullet X$.
\item
When $\cat{C}=\cat{0}$ (the initial category) and $\cat{C'}=\cat{C''}=\cat{1}$, then $\Psi$ is a presheaf $X$, $\Phi$ is a set $z$, and the definition boils down to ${\it Lan}_\subseteq X z$.
\end{enumerate}
It is also easily checked that the identity $I$ as defined in Section \ref{Kelly-op} is indeed the identity of $\cat{Prof}_?$.

\medskip
We end the section by giving a formula for $\Psi\bullet\Phi$ when $!$ is now an {\em arbitrary} monad $(!,\eta_!,\mu_!)$ over \cat{Cat} given together with a distributive law: $\dist:\: !\comp{\it Psh}\rightarrow\:{\it Psh}\:\comp\: !$. This could open the way for 
handling variations on the shapes of operations (the third kind of variation considered in the introduction) in a profunctor setting. 

Recall (Section \ref{Profunctors}) that given $\Phi:\cat{C}\profarrow\cat{C'}$ we write
$\Phi'$ for its presentation as a functor from $\cat{C'}^{\it op}$ to $\cat{Set}^{\cat{C}}$.
Then we can give abstract versions of the steps $(1)^{\it op}$, $(2)^{\it op}$ and $(3)$ above, as follows:
(we still write $?\cat{C}=(!\cat{C}^{\it op})^{\it op}$):

\medskip\noindent
$(1)^{\it op}$ Given $\Phi:\cat{C}\profarrow\cat{C'}$, we can define $?\Phi:?\cat{C}\profarrow ?\cat{C'}$ by the equation:
$$(?\Phi)'=\dist\comp !(\Phi')$$
where on the right hand side we apply $!$ as a functor on \cat{Cat}.

\medskip\noindent
$(2)^{\it op}$ The comultiplication $\delta:?\cat{C}\profarrow ??\cat{C}$ is defined by the formula
$$\delta'=\eta_{{\it Psh}}\comp (\mu_!)_{\it \cat{C}^{\it op}}$$
 
 \medskip\noindent
(3) The composition in \cat{Prof} of  $\Phi':\cat{C'}^{\it op}\rightarrow\cat{Set}^{\cat{C}}$ and
$\Psi':\cat{C''}^{\it op}\rightarrow\cat{Set}^{\cat{C'}}$ is defined as $\Phi'^\#\comp\Psi'$.
Hence, given $\Phi:?\cat{C}\profarrow \cat{C'}$, we obtain the following formula for $\Phi^{\flat}:?\cat{C}\profarrow ?\cat{C'}$:
$${\Phi^{\flat}}'=(\eta_{{\it Psh}}\comp (\mu_!)_{\it \cat{C}^{\it op}})^{\#}\comp\dist\comp(!\Phi')=
{\it Psh}( (\mu_!)_{\it \cat{C}^{\it op}})\comp\dist\comp(!\Phi')$$

\medskip\noindent
Finally, expanding the definition of ${\it Psh}$, we get the following formula for the composition $\Psi\bullet\Phi=\Psi\comp\Phi^{\flat}$ in
$\cat{Prof}_?$:

$$(\Psi\bullet\Phi)(\gamma,C'')=\int^{\gamma'}\Psi(\gamma',C'')\times\int^{\delta}\dist((!\Phi')\gamma')\delta\times ?\cat{C}[\mu_!\delta,\gamma]$$

\medskip
Following the style of~[\cite{Tanaka-Power}]~, one may also hide the distributive law and write

$$(\Psi\bullet\Phi)(\gamma,C'')=\int^{\gamma'}\Psi(\gamma',C'')\times\overline{\gamma'}(\Phi')\gamma$$
where $\overline{\gamma'}(\Phi')\gamma=\hat{\Phi}(\gamma,\gamma')$ can be described abstractly as follows.  For each (pseudo-)$!$-algebra $(\cat{D},\alpha:!\cat{D}\rightarrow \cat{D})$, for each category $\cat{C}$, and each object $C$ of $!\cat{C}$, $C$ induces a functor $\overline{C}=\cat{D}^\cat{C}\rightarrow\cat{D}$ defined as $(F\mapsto \alpha(!FC))$. The above formula is obtained by instantiating $\cat{C}$, $C$, $\cat{D}$, and $\alpha$ as
$\cat{C'}^{\it op}$, $\gamma'$, $\cat{Set}^{!(\cat{C}^{\it op})}$, and ${\it Psh}(\mu_!)\comp\dist$, respectively. That this is a (pseudo)-$!$-algebra structure is established using the equalities satisfied by $\dist$.

\medskip
This leads us to a generalised definition of operad.

\begin{definition}
Given a monad $!$ on \cat{Cat} and a distributive law $\dist:\: !\comp{\it Psh}\rightarrow\:{\it Psh}\:\comp\: !$, a $\cat{C}$-coloured $!$-operad is a monoid in $\cat{Prof}_?[\cat{C},\cat{C}]$, that is, a functor $X:\:?\cat{C}\times\cat{C}^{\it op}\rightarrow\cat{Set}$ given together with two natural transformations $e:{\it id}\rightarrow X$ and $m:X\bullet X\rightarrow X$ satisfying the monoid laws.  A $\cat{1}$-coloured $!$-operad is called a $!$-operad.
\end{definition}

Hence the non-symmetric operads (resp. symmetric operads, clones) are the $\cat{1}$-coloured $!_\emptyset$-operads
(resp. $!_{\set{\sigma}}$-operads, $!_{\set{\sigma,\delta,\epsilon}}$-operads), and
the non-symmetric coloured operads (resp. symmetric coloured operads) are the coloured $!_\emptyset$-operads
(resp. coloured $!_{\set{\sigma}}$-operads). Indeed, coloured operads were defined in this way by Dolan and Baez~[\cite{BD98}]~(though they do not explicate the underlying distributive law or lifting, cf. Footnote \ref{BD-not-dist}).

\section{Cooperads and properads}

So far, we have addressed the variation on the ``first axis" that goes from non-symmetric operads to clones.
In this section, we address the variation on the second axis (cooperations, bioperations).

Cooperations are dual to operations. Sets of cooperations are organised into cooperads, the notion dual to that of operad. One works now in the Kleisli category $\cat{Prof}_!$, i.e.,
a cooperad is a profunctor $X:\cat{1}\profarrow!\cat{1}$ with a monoid structure in $\cat{Prof}_!,$, and we can vary on the first axis as we did for operads.

For bioperations, the natural idea is to consider profunctors from $?\cat{1}$ to $!\cat{1}$. Such profunctors 
can compose provided there is a distribitutive law (another one!) $\dist: \:?_X!_X\rightarrow \:!_X?_X$: given $\Phi,\Psi:?\cat{1}\profarrow !\cat{1}$, we define $\Psi\bullet\Phi$ 
by composing the comultiplication of $?$, $?\Phi$, $\dist$, $!\Psi$, and the multiplication of $!$. 
The identity is the composition ot the counit of $?$ and of the unit of $!$.

This idea has been carried out in detail
by Garner, in the case of $X=\set{\sigma}$~[\cite{Garner}]~.  We set 
$$\dist((m_1,\ldots,m_p),(n_1,\ldots,n_q))\neq\emptyset\quad\mbox{iff}\quad m_1+\ldots+m_p=
n_1+\ldots+n_q$$
and when the equality holds, $\dist((m_1,\ldots,m_p),(n_1,\ldots,n_q))$ is the set of 
permutations $s$ from  $m_1+\ldots+m_p$ to $n_1+\ldots+n_q$ such that the graph 
\begin{itemize}
\item whose 
set of vertices is the disjoint union of
$$\set{1,\ldots,p}, \set{1,\ldots m_1+\ldots+m_p}, \set{1,\ldots,n_1+\ldots+n_q},\mbox{ and }\set{1,\ldots,q}$$
 \item and whose set of edges is the union 
 \begin{itemize}
 \item of the edges between $1$ and $1,\ldots m_1$, between $2$ and $m_1+1,\ldots,m_2$, \ldots, and between 
 $p$ and $m_{p-1}+1,\ldots,m_p$,
 \item of the edges between $i$ and $s(i)$ ($i\leq m_1+\ldots+m_p$), and
 \item of the edges between $1,\ldots n_1$ and $1$, \ldots, and between 
 $n_{q-1}+1,\ldots,n_q$ and $q$.
 \end{itemize}
 \end{itemize}
 is connected.  This formalises the idea that when composing two sets of bioperations, we are interested primarily in the compositions that preserve connectedness. More precisely, composition of properations has both a vertical aspect and a horizontal aspect. The distributive law takes care of the vertical aspect. 
Disjoint connected compositions can be placed in parallel and composed horizontally. Thus, a full categorical account requires a double category setting, that takes care of these two aspects. We refer to~[\cite{Garner}]~for details.

Monoids for this composition operation provide an analogue of operads for bioperations, and for $X=\set{\sigma}$ and $\dist$ as described above, what we obtain (replacing \cat{Set} by \cat{Vect}) is exactly the notion of {\em properad} introduced by Vallette in his Th\`ese de Doctorat~[\cite{ValletteTh, ValletteAMS}]~.

\paragraph{Acknowledgements.}
I collected the material presented  here for an invited talk at the conference Operads 2006.
When preparing this talk, I benefited a lot
from discussions with Marcelo Fiore and Martin Hyland. The string diagrams have been drawn with {\tt strid}, a tool due to Samuel Mimram and Nicolas Tabareau (\url{http://strid.sourceforge.net/}). I also wish to thank the anonymous referee for his helpful remarks. He in particular pointed out an alternative way to address the size issue, namely to restrict attention to presheaves that are limits of small diagrams of representables.

\end{document}